\documentclass[10pt,oneside]{amsart}
\usepackage{latexsym,amssymb,amsfonts,amsmath,graphicx}

\makeatletter
\renewcommand{\@secnumfont}{\bfseries}
\makeatother

\newcommand{\QED}{\hfill $\Box$}
\newcommand{\N}{\mathbb N}
\newcommand{\Z}{\mathbb Z}
\newcommand{\C}{\mathbb C}
\newcommand{\un}[1]{\underline{#1}}

\begin{document}
\thispagestyle{empty}

\begin{center}
{\Large \bf 
A structured description of the \\ \vspace*{0.3em}
genus spectrum of abelian $p$-groups}
\end{center}

\vspace*{1em}

\begin{center}
{\large\bf J\"{u}rgen M\"{u}ller and Siddhartha Sarkar}
\end{center}

\vspace*{1em}

{\small \begin{center}{\bf Abstract}\end{center}

The genus spectrum of a finite group $G$ is the set of all $g$ such that
$G$ acts faithfully on a compact Riemann surface of genus $g$. It is an
open problem to find a general description of the genus spectrum of the
groups in interesting classes, such as the abelian $p$-groups. Motivated
by the work of Talu \cite{tal} for odd primes $p$, we develop a general
combinatorial machinery, for arbitrary primes, to obtain a structured
description of the so-called reduced genus spectrum of abelian $p$-groups.

We have a particular view towards how to generally find the reduced 
minimum genus in this class of groups, determine the complete genus
spectrum for a large subclass of abelian $p$-groups, consisting of
those groups in a certain sense having `large' defining invariants,
and use this to construct infinitely many counterexamples to 
Talu's Conjecture \cite{tal}, saying that an abelian $p$-group is
recoverable from its genus spectrum. Finally, we indicate the
effectiveness of our combinatorial approach by applying it to 
some explicit examples.

Mathematics Subject Classification (MSC2010): 
20H10; 20K01, 30F35, 57S25. 
}

\begin{center}\small {\bf Contents} 

\begin{tabular}{rp{240pt}r}
\ref{introsec} & Introduction \dotfill & \pageref{introsec}\\
\ref{presec} & Groups acting on Riemann surfaces \dotfill & \pageref{presec}\\
\ref{mainlinesec} & Mainline integers \dotfill & \pageref{mainlinesec}\\
\ref{talusec} & Talu's Theorem revisited \dotfill & \pageref{talusec}\\
\ref{transsec} & Transforming to mainline integers \dotfill 
               & \pageref{transsec}\\
\ref{proofsec} & The main result \dotfill & \pageref{proofsec}\\
\ref{applsec} & Talu's Conjecture \dotfill & \pageref{applsec}\\
\ref{smallranksec} & Examples: Small rank \dotfill & \pageref{smallranksec}\\
\ref{smallexpsec} & Examples: Small exponents \dotfill & \pageref{smallexpsec}\\
\end{tabular}
\end{center}


\section{\bf Introduction}\label{introsec}

\subsection{Genus spectra}
Given a compact Riemann surface $X$ of genus $g\geq 0$, 
a finite group $G$ is said to act on $X$, if $G$ can be embedded
into the group $\mathrm{Aut}(X)$ of biholomorphic maps on $X$. 
While $\mathrm{Aut}(X)$ is infinite as long as $g\leq 1$,
by the Hurwitz Theorem \cite{hur} we have 
$|\mathrm{Aut}(X)|\leq 84\cdot(g-1)$ as soon as $g\geq 2$. 
Thus in the latter case there are only finitely many groups $G$,
up to isomorphism, acting on $X$.

But conversely, given a finite group $G$ there always is an infinite set
$\mathrm{sp}(G)$ of integers $g\geq 0$, called the {\bf (genus) spectrum}
of $G$, such that there is a Riemann surface $X$ of genus $g$
being acted on by $G$; in this case, $g$ is called a genus of $G$.
Note that we are in particular including the cases $g\leq 1$.
In \cite{mmi}, the problem of determining $\mathrm{sp}(G)$ is called
the {\bf Hurwitz problem} associated with $G$, and the problem of 
finding the {\bf minimum genus} $\min~\mathrm{sp}(G)$ of $G$, also called its
{\bf strong symmetric genus}, has arisen some particular interest. 
For more details
we refer the reader to \cite{bre,sar}, and the references given there.

To attack the Hurwitz problem, let $\Delta(G):=\frac{|G|}{\mathrm{exp}(G)}$,
where $\mathrm{exp}(G)$ denotes the exponent of $G$,
that is the least common multiple of the orders of its elements.
Then let the {\bf reduced (genus) spectrum} of $G$ be defined by
$$ \mathrm{sp}_0(G):=\left\{ \frac{g-1}{\Delta(G)}\in\Z ~:~
g\in\mathrm{sp}(G)\right\} ,$$
where the number $\frac{g-1}{\Delta(G)}$ 
is called the reduced genus associated with $g$.
It follows from \cite{kul}, together with a special
consideration of the case $g=0$, that
$$ \mathrm{sp}_0(G) \subseteq \mathbb{S}:=
\frac{1}{\epsilon(G)}\cdot(\{-1\}\cup\N_0) $$
is a co-finite subset, where $\epsilon(G)$ divides $\gcd(2,|G|)$ 
and can be determined from the structure of $G$, 
as is recalled in (\ref{kulkarnithm}).
A word of caution is in order here: In \cite{kul} the
notion of reduced genus is defined differently, by 
taking $\epsilon(G)$ into account as well, while 
our choice leads to fewer case distinctions.

The {\bf reduced minimum genus} of $G$, that is the reduced
genus associated with the minimum genus of $G$, equals 
$\mu_0(G):=\min~\mathrm{sp}_0(G)$. 
Moreover, following \cite{kma}, the {\bf reduced stable upper genus} 
$\sigma_0(G)$ of $G$ is the smallest element of $\mathbb{S}$ such that
all elements of $\mathbb{S}\setminus\mathrm{sp}_0(G)$ are less than
$\sigma_0(G)$; the genus $\sigma(G)$ associated with $\sigma_0(G)$
is called the {\bf stable upper genus} of $G$.
The elements of $\mathbb{S}\setminus\mathrm{sp}_0(G)$ exceeding $\mu_0(G)$ 
are called the {\bf reduced spectral gap}
of $G$; the associated genera form the {\bf spectral gap} of $G$.
Hence solving the Hurwitz problem for $G$ amounts to determining
$\mu_0(G)$ and $\sigma_0(G)$ and the reduced spectral gap of $G$.


\subsection{Our approach to abelian $p$-groups}
We now restrict ourselves to finite $p$-groups $G$,
where $p$ is a prime. Not too much is known about
the genus spectrum of groups within this class, not even if we only
look at interesting subclasses, for example those given by bounding
a certain invariant such as rank, exponent, nilpotency class,
or co-class; see \cite{sar}. 

This still holds if we restrict further to the class of abelian $p$-groups,
which are the groups we are interested in from now on, 
their general shape being
\[ G \cong {\mathbb Z}_p^{r_1} \oplus {\mathbb Z}_{p^2}^{r_2} 
  \oplus \cdots \oplus {\mathbb Z}_{p^e}^{r_e}, \]
where $e\geq 1$, and $r_i\geq 0$ for $1\leq i\leq e-1$, and $r_e \geq 1$.
We point out that, in particular contrary to \cite{mta,tal},
we are allowing for arbitrary primes $p\geq 2$ throughout.

We give an outline of the paper:
In {\bf Section \ref{presec}} we recall a few facts 
about Riemann surfaces and their automorphism groups.
In {\bf Section \ref{mainlinesec}} we prepare the combinatorial
tools needed later on; we comment on them in (\ref{mainlinecomment}).
Having these preliminaries in place we turn out attention to
abelian $p$-groups and their genus spectra:

{\bf Section \ref{talusec}:}
Our starting point is Talu's approach \cite{tal}
towards a general description of the genus spectrum of abelian $p$-groups,
in the case where $p$ is odd. Building on these ideas, we develop a 
conceptual approach to describe the smooth epimorphisms, in the sense 
of (\ref{fuchs}), onto a given abelian $p$-group, where $p$ is arbitrary.
The resulting general necessary and sufficient arithmetic condition for 
their existence, which we still refer to as {\bf Talu's Theorem}, 
is given in Theorems (\ref{smoothnecthm}) and (\ref{smoothsuffthm}); 
in proving the latter we in particular close a gap 
in the proof of \cite[Thm.3.3]{tal}.

{\bf Section \ref{transsec}:}
This is then translated into a
combinatorial description of the domain of the reduced genus map,
yielding a structured description of the reduced spectrum of $G$ 
being presented in (\ref{talurem}),
and leading to a machinery to compute the 
reduced minimum genus $\mu_0(G)$ of $G$ culminating in 
Theorem (\ref{minthm}), which says that $\mu_0(G)$ is given as 
the minimum of at most $e+1$ numbers, given explicitly
in terms of the defining invariants $(r_1,\dotsc,r_e)$.
In particular, in (\ref{macrem}) we obtain
an independent proof and an improved version of Maclachlan's
method \cite[Thm.4]{mac} for the special case of abelian $p$-groups. 
Our combinatorial approach should also be suitable to get 
hands on the reduced stable upper genus $\sigma_0(G)$ of $G$;
we are planning to pursue this further in a subsequent paper.

{\bf Section \ref{proofsec}:}
Having this combinatorial machinery in place, we turn to abelian
$p$-groups with `large' invariants, by assuming that
$$ r_i\geq p-1 \quad \textrm{for} \quad 1\leq i\leq e-1, 
\quad \textrm{and} \quad r_e \geq \max\{p-2,1\} .$$
In these cases we are able to determine both the reduced minimum 
genus $\mu_0(G)$ as well as the reduced stable upper genus $\sigma_0(G)$ 
in terms of the defining invariants $(r_1,\dotsc,r_e)$ of $G$. 
More precisely, our main result says the following:

\subsection*{Main Theorem (\ref{genusthm})}
Let $G$ have `large' invariants as specified above.
Then the reduced minimum and stable upper genera of $G$ are given as
\[ \mu_0(G)= \sigma_0(G) =
\frac{1}{2}\cdot\bigg(-1-p^e+\sum_{i=1}^e (p^e-p^{e-i})\cdot r_i \bigg) .\]


At this stage, a comparison with \cite{tal} is in order: 
The major aim there is to study abelian $p$-groups
having `small' invariants, fulfilling
$1+\sum_{j=i}^er_j\leq (e-i+1)\cdot (p-1)$, for $1\leq i\leq e$,
with a particular view towards computing the reduced stable upper
genus $\sigma_0(G)$ in these cases, the key result being a 
closed formula for $\sigma_0(G)$ in terms of the defining invariants 
$(r_1,\dotsc,r_e)$.
Now one of the maximal admissible `small' cases 
coincides with the smallest admissible case here, 
thus we recover \cite[Cor.3.7]{tal}, where $\sigma_0(G)$
is explicitly determined, but $\mu_0(G)$ is only claimed without proof.

{\bf Section \ref{applsec}:}
Next, we turn to an aspect of the general question of how much
information about a group is encoded into its spectrum, at best
whether its isomorphism type can be recovered from it.
Since in view of the examples in \cite{mta} this cannot possibly 
hold without restricting the class of groups considered, the class
of abelian $p$-group seems to be a good candidate to look at. 
More specifically, {\bf Talu's Conjecture} \cite{tal} says
that, whenever $p$ is odd, the spectrum of a non-trivial abelian $p$-group 
already determines the group up to isomorphism. Moreover,
although this cannot possibly hold in full generality for $p=2$, 
for example in view of the sets of groups $\{\Z_2,\Z_4,\Z_2^2,\Z_8\}$ and
$\{\Z_2\oplus\Z_4,\Z_2^3,\Z_2\oplus\Z_8\}$ discussed below, 
we are tempted to expect that it still holds true up to 
finitely many finite sets of exceptions.  

But, as a consequence of (\ref{genusthm}), we are able to produce
infinitely many counterexamples to Talu's Conjecture (both for $p$ odd
and $p=2$), that is pairs of non-isomorphic abelian
$p$-groups having the same spectrum. We present two distinct kinds of
counterexamples, consisting of groups having the same order and exponent,
and of groups where these invariants are different, in (\ref{e3ex})
and (\ref{counterex}), respectively.
This also shows that there cannot be an absolute bound on the cardinality
of a set of abelian $p$-groups sharing one and the same spectrum, not even
if we restrict ourselves to groups having the same order and exponent.
Still, we will have to say something positive on Talu's Conjecture later on.

{\bf Section \ref{smallranksec}:}
In order to show the effectiveness of the combinatorial machinery developed
we work out various examples, where in particular we get new systematic
proof of a number of earlier results scattered throughout the literature:
In (\ref{nonposex}) we determine the groups of non-positive reduced
minimum genus, where we recover the abelian $p$-groups amongst the
well-known finite groups acting on surfaces of genus $g\leq 1$, 
see \cite[App.]{sah} or \cite[Sect.6.7]{dou}. In particular, the 
non-cyclic abelian groups of order at most $9$, which have to be
treated as exceptions in \cite[Thm.4]{mac}, reappear here naturally.

In (\ref{rank2}) we deal with the groups of rank at most $2$, 
whose smallest positive reduced genus we determine. 
In particular, for the cyclic groups we recover the results in 
\cite{har} and \cite[Prop.3.3]{kul}, for the groups of rank $2$ 
we improve the bound in \cite[Prop.3.4]{kul}, 
and for the cases of cyclic deficiency $1$, where $p$ is odd, we recover
the relevant part of \cite[Thm.5.4]{mta} and \cite[Cor.5.5]{mta}.
Moreover, we show that a cyclic $p$-group is uniquely
determined by its smallest genus $\geq 2$, with the single
exception of the groups $\{\Z_2,\Z_4,\Z_8\}$, and 
that an abelian $p$-groups of rank $2$ 
is uniquely determined by its smallest genus $\geq 2$, with the 
single exception of the groups $\{\Z_2\oplus\Z_4,\Z_2\oplus\Z_8,\Z_4^2\}$.

{\bf Section \ref{smallexpsec}:}
In (\ref{e1ex}) and (\ref{e2ex}) we determine the reduced minimum genus
of the elementary abelian $p$-groups, and of the abelian $p$-groups of
exponent $p^2$, respectively. Using this, we show that within the class 
of elementary abelian $p$-groups a group is uniquely determined by its 
minimum genus, with the single exception of the groups $\{\Z_2,\Z_2^2\}$;
for $p$ odd this would also be a consequence of \cite[Cor.7.3]{mta}, 
but \cite[Sect.7, Rem.]{mta} preceding it contains an error.
Similarly, we show that within the class of abelian $p$-groups of
exponent $p^2$ a group is uniquely determined by its Kulkarni invariant,
see (\ref{kulkarnithm}), and its minimum genus,
with the single exception of the groups $\{\Z_4^2,\Z_2\oplus\Z_4\}$;  
for $p$ odd this is claimed without proof in \cite[Thm.3.8]{tal}.

To summarize our results in Sections \ref{smallranksec} and \ref{smallexpsec}, 
although Talu's Conjecture is false in general, it turns out to hold 
within the following subclasses of the class of non-trivial
abelian $p$-groups (including the case $p=2$):
{\bf i)} the class of cyclic $p$-groups,
{\bf ii)} the class of $p$-groups of rank $2$,
{\bf iii)} the class of elementary abelian $p$-groups, and
{\bf iv)} the class of $p$-groups of exponent $p^2$.

\subsection{Mainline integers}\label{mainlinecomment}
We comment on the combinatorial tool featuring prominently
in our approach: Given a prime $p$, and a non-increasing 
sequence $\un{a}:=(a_1,\dotsc,a_e)$ of non-negative integers,
the associated {\bf $p$-mainline integer} 
(as we call it by lack of a better name) 
is defined as $\wp(\un{a}):=\sum_{i=1}^e a_ip^{e-i}$.
Moreover, given any non-increasing sequence $\un{s}:=(s_1\dotsc,s_e)$ 
of non-negative integers, let ${\mathcal P}(\un{s})$ be the set of all
$p$-mainline integers $\wp(\un{a})$ where $\un{a}$ is bounded below 
component-wise by $\un{s}$.  
The connection to abelian $p$-groups with defining
invariants $(r_1,\dotsc,r_e)$ is given by letting the sequence $\un{s}$
be given by
$$ s_i:=1+\sum_{j=i}^e r_j\quad\textrm{for}\quad 1\leq i\leq e .$$

We are interested in the structure of ${\mathcal P}(\un{s})$, whose
minimum obviously equals $\wp(\un{s})$. It can be shown that 
${\mathcal P}(\un{s})$ is a co-finite subset of the non-negative integers, 
and thus the combinatorial problems arising are to 
determine the smallest $m$ such that all integers from $m$ on actually
are elements of ${\mathcal P}(\un{s})$, and to describe the gap 
consisting of the non-mainline integers between $\wp(\un{s})$ and $m$.

It might very well be possible that this general kind of  
problems is well-known to combinatorialists, but we have not been
able to find suitable references. In consequence we develop
a piece of theory, just as far as necessary for the present paper;
we are planning to elaborate on this, as we go along with
pursuing further questions concerning the genus spectrum of abelian
$p$-groups.

\section{\bf Groups acting on Riemann surfaces}\label{presec}

We assume the reader familiar with the basic theory of Riemann surfaces,
as is exhibited for example in \cite{bre,dou}, so that here we are 
just content with recalling a few facts.
The connection between geometry and group theory is given by
the following well-known theorem. We point out that it is often 
only used for $g\geq 2$, in which case the `groups with signature'
occurring are the Fuchsian groups, but it actually holds for 
all $g\geq 0$; see for example 
\cite[Sect.1]{bre} and \cite[Ch.6]{dou} and \cite{sah}:

\subsection{Theorem}\label{fuchsthm}
A finite group $G$ acts on a compact Riemann surface $X$,
if and only if there is $\Gamma\leq\mathrm{Aut}(U)$, 
where $U$ is a simply-connected Riemann surface and
$\Gamma$ is a group with signature in the sense of (\ref{fuchs}), 
and a smooth epimorphism $\phi : \Gamma \longrightarrow G$,
such that $X$ is isomorphic to the orbit space $U/\mathrm{ker}(\phi)$.
\QED

\subsection{Smooth epimorphisms}\label{fuchs}
We keep the notation of (\ref{fuchsthm}).
A group $\Gamma$ is said to be a {\bf group with (finite) signature} 
if it has a distinguished generating set
$$ \{a_k,b_k ~:~ 1\leq k\leq h\}\quad\cup\quad\{c_j ~:~ 1\leq j\leq s\} ,$$
for some $h,s\in\N_0$, subject to the order relations
$$  c_j^{n_j} = 1, \quad\textrm{where}\quad n_j\in\N\setminus\{1\}, $$
for $1\leq j\leq s$, and the `long' relation
$$ \prod_{k=1}^{h} [a_k, b_k]\cdot\prod_{j=1}^s c_j = 1 ,$$
where $[a,b]:=a^{-1}b^{-1}ab$ denotes the commutator of $a$ and $b$.
More generally, there might also be order relations of the form
`$c^\infty=1$', that is no order relation for the generator $c$ at all; 
but since we are requiring $X$ to be compact, and hence the orbit space 
$X/G$ to be compact as well, these cases do not occur here;
see \cite[App.]{sah}.

An epimorphism $\phi : \Gamma \longrightarrow G$
with torsion-free kernel is called {\bf smooth}.
This is equivalent to the condition that
$$ \phi(c_j)\in G\quad\textrm{has order}\quad n_j,
\quad\textrm{for all}\quad 1\leq j\leq s .$$ 
In this case, the $(s+1)$-tuple $(n_1,\dotsc,n_s;h)$ is called a
{\bf signature} of $G$, with {\bf periods} $n_1,\dotsc,n_s\geq 2$ 
and {\bf orbit genus} $h\geq 0$.
The orbit space $X/G$ has genus $h$, and the branched covering
$X\longrightarrow X/G$ gives rise to the {\bf Riemann-Hurwitz equation}
$$ g-1 =|G|\cdot\bigg(h-1
+\frac{1}{2}\cdot\sum_{i=1}^s(1-\frac{1}{n_i})\bigg) .$$

\subsection{Kulkarni's Theorem}\label{kulkarnithm}
To describe the structure of the genus spectrum of a finite group $G$,
in \cite{kul} a group theoretic invariant $N(G)\in\N$, now called the 
{\bf Kulkarni invariant} of $G$, is introduced, such that
$$ \mathrm{sp}(G)\setminus\{0\} \subseteq 1+ N(G)\cdot\N_0 ,$$
and $\mathrm{sp}(G)\setminus\{0\}$ is a co-finite subset of 
$1+N(G)\cdot\N_0$. 
Moreover, we have
$$ N(G) = \frac{1}{\epsilon(G)}\cdot\frac{|G|}{\mathrm{exp}(G)} ,$$
where $\epsilon=\epsilon(G)\in\{1,2\}$ is determined 
by the structure of $G$ as follows:

If $|G|$ is odd, then $\epsilon:=1$; if $|G|$ is even, letting $\tilde G$ 
be a Sylow $2$-subgroup of $G$, then $\epsilon:=1$ provided the subset
$\{a\in \tilde G;|a|<\mathrm{exp}(\tilde G)\} \subseteq G$
forms a subgroup of $\tilde G$ of index $2$, otherwise $\epsilon:=2$.
In other words, using the notions developed in {\cite{ms}}, we have 
$\epsilon=2$ if and only if $\tilde G$ is a non-trivial $2$-group 
not of `GK type'. 

This yields the description of the non-negative part
of the reduced spectrum $\mathrm{sp}_0(G)$ as stated earlier.
As for its negative part, the well-known description of finite group
actions on compact Riemann surfaces of genus $g=0$,
see \cite[App.]{sah} or \cite[Sect.6.7]{dou}, 
says that in this case $G$ is cyclic, dihedral, alternating or symmetric
of isomorphism type in
$\{\Z_n,\textrm{Dih}_{2n},\textrm{Alt}_4,\textrm{Sym}_4,\textrm{Alt}_5\}$,
hence we indeed get $\Delta(G)=\epsilon(G)$.

\subsection{The case of $p$-groups}\label{genusmap}
We turn to the case of interest for us:
Let $G$ be a $p$-group of order $p^n$ and exponent $p^e$, 
where $e\leq n\in\N_0$.

If $\phi : \Gamma \longrightarrow G$ is a smooth epimorphism,
then all the periods are of the form $p^i$, where 
$0\leq i\leq e$. Hence we may abbreviate any
signature $(n_1,\dotsc,n_s;h)$ of $G$ by the $(e+1)$-tuple
$(x_1,\dotsc,x_e;h)$, being called the associated {\bf $p$-datum}, where 
$$ x_i:=|\{1\leq j\leq s; n_j=p^i\}|\in\N_0 .$$
 
The set $D(G)$ of all $p$-data of $G$, being afforded by smooth 
epimorphisms, is called the {\bf data spectrum} of $G$.
Then the Riemann-Hurwitz equation gives rise to the 
{\bf genus map} $g : D(G) \longrightarrow \mathrm{sp}(G)$ defined by
$$ g(x_1,\dotsc,x_e;h) :=
1+ p^{n-e} \cdot \bigg((h-1)\cdot p^e + \frac {1}{2}\cdot \sum_{i=1}^e 
x_i ( p^e - p^{e-i}) \bigg) . $$ 

Letting the {\bf cyclic deficiency} of $G$ be defined as 
$$ \delta=\delta(G):=\log_p(\Delta(G))=n-e\in\N_0 ,$$
in view of Kulkarni's Theorem (\ref{kulkarnithm}) we have 
$N(G) = \frac{1}{\epsilon(G)}\cdot p^{\delta(G)}$. 
Then the {\bf reduced genus map} 
$g_0 : D(G) \longrightarrow \mathrm{sp}_0(G)
\subseteq \frac{1}{\epsilon(G)}\cdot(\{-1\}\cup\N_0)
\subseteq \frac{1}{2}\cdot(\{-1\}\cup\N_0)$,
given by associating the reduced genus
$\frac{g-1}{p^\delta}\in\mathrm{sp}_0(G)$ 
with any $g\in \mathrm{sp}(G)$, reads
\[ g_0(x_1, \dotsc, x_e;h) = 
(h-1)\cdot p^e + {\frac {1}{2}} \cdot\sum_{i=1}^{e} x_i(p^e - p^{e-i}) .\]

\section{\bf Mainline integers}\label{mainlinesec}

In this section we consider sequences of non-negative integers
from a certain purely combinatorial viewpoint. We develop a 
little piece of general theory, as far as will be needed in 
Sections \ref{transsec} and \ref{proofsec}.

\subsection{Integer sequences}\label{mainlinesetting}
Given finite sequences $\un{a}=(a_1,\dotsc,a_e)\in{\mathbb N}_0^e$
and $\un{b}=(b_1,\dotsc,b_e)\in{\mathbb N}_0^e$ of non-negative integers,
of length $e\geq 1$, we write $\un{a}\leq\un{b}$, 
and say that $\un{b}$ {\bf dominates} $\un{a}$,
if $a_i \leq b_i$ for all $1\leq i\leq e$.
We will be mainly concerned with the set of {\bf non-increasing} sequences
$$ {\mathcal N}={\mathcal N}(e):=\{\un{a}=(a_1,\dotsc,a_e)\in{\mathbb N}_0^e
~:~ a_1\geq \cdots \geq a_e \} .$$
We introduce a few combinatorial notions concerning integer sequences:
To this end, we fix $p\in\N$; later on $p$ will be a prime, but here
is no need to assume this.

{\bf i)}
For an arbitrary sequence $\un{a}=(a_1,\dotsc,a_e)\in {\mathbb N}_0^e$ let
$$ \wp(\un{a})=\wp(a_1,\dotsc, a_e):= 
\sum_{i=1}^e a_ip^{e-i} \in {\mathbb N}_0 .$$
Then the {\bf ($p$-)mainline integers} associated with $\un{a}$ are defined as 
$$ {\mathcal P}(\un{a})={\mathcal P}(a_1,\dotsc,a_e) := 
\{ \wp(\un{b})\in {\mathbb N}_0 ~:~ 
\un{b}\in{\mathcal N},\,\un{a}\leq\un{b} \}. $$
Note that we allow for arbitrary $\un{a}$ to start with,
while the sequences $\un{b}$ used in the definition of ${\mathcal P}(\un{a})$
are required to be non-increasing. It will turn out that there
always is a non-increasing sequence affording a given set of
mainline integers.

The {\bf hull sequence} 
$\un{\tilde a}=(\tilde a_1,\dotsc,\tilde a_e)\in {\mathcal N}$ 
of $\un{a}$ is defined recursively by letting $\tilde a_e:=a_e$ and
$$ \tilde a_i:=\mathrm {max}\{\tilde a_{i+1},a_i\}\quad \textrm{for} \quad
e-1\geq i\geq 1 ;$$
note that this definition is actually independent of the chosen integer $p$.
Hence we have $\un{a}\leq\un{\tilde a}$, where 
$\un{a}=\un{\tilde a}$ if and only if $\un{a}\in{\mathcal N}$.

{\bf ii)}
Given a non-increasing sequence $\un{a}=(a_1,\dotsc,a_e)\in{\mathcal N}$,
its {\bf $p$-enveloping sequence}
$\un{\hat a}=(\hat a_1,\dotsc,\hat a_e)\in{\mathcal N}$ 
is defined recursively by $\hat a_e:=a_e$ and 
$$ \hat a_i := {\mathrm {max}} \{ \hat a_{i+1} + (p-1), a_i \} 
\quad \textrm{for} \quad e-1 \geq i \geq 1 ;$$
hence we have $\un{a}=\un{\tilde a}\leq\un{\hat a}$,
where $\un{a}=\un{\hat a}$ if $p=1$.

Moreover, whenever $e\geq 2$ let
$$ ||\un{a}|| = ||(a_1,\dotsc, a_e)|| := 
{\mathrm {min}} \{ a_i - a_{i+1} ~:~ 1 \leq i \leq e-1 \} ,$$
and let $||\un{a}||:=\infty$ for $e=1$; note that, despite notation,
$||\cdot ||$ is not a norm in sense of metric spaces. In particular, 
we have $\un{a}=\un{\hat a}$ if and only if $||\un{a}||\geq p-1$.

\subsection{Proposition}\label{mainlineprop}
Given $\un{a}\in{\mathbb N}_0^e$,
then we have ${\mathcal P}(\un{a})={\mathcal P}(\un{\tilde a})$.

\subsection*{Proof}
Let $\un{b}=(b_1,\dotsc, b_e)\in{\mathcal N}$.
If $\un{\tilde a}\leq\un{b}$, then from $\un{a}\leq\un{\tilde a}$
we also get $\un{a}\leq\un{b}$. Conversely, if $\un{a}\leq\un{b}$,
then we have $\tilde a_e=a_e\leq b_e$, and recursively for
$e-1\geq i\geq 1$ we get
$\tilde a_{i+1}\leq b_{i+1}\leq b_i$ and $a_i\leq b_i$, hence
$\tilde a_i\leq b_i$; this implies that $\un{\tilde a}\leq\un{b}$.
\QED

\subsection{Proposition}\label{cofinprop}
Given $\un{a}\in{\mathcal N}$, the set
${\mathbb N}_0\setminus {\mathcal P}(\un{a})$ is finite.

\subsection*{Proof} 
We consider the $p$-enveloping sequence 
$\un{\hat a}=(\hat a_1,\dotsc,\hat a_e)\in{\mathcal N}$ 
of $\un{a}$, and we show that any $m\geq\wp(\un{\hat a})$ is a 
mainline integer: To this end, write $m-\wp(\un{\hat a})$ in 
a partial $p$-adic expansion as
$m-\wp(\un{\hat a}) = \sum_{i=1}^e b_ip^{e-i}$,
where $b_i\geq 0$ such that $b_2, \dotsc, b_e \leq p-1$, 
but $b_1$ might be arbitrarily large.
Then we have
$m = \sum_{i=1}^e (\hat a_i+b_i)p^{e-i}$.
Since for $1\leq i\leq e-1$ we have
$\hat a_i - \hat a_{i+1} \geq p-1 \geq b_{i+1} - b_i$, thus
$\hat a_i + b_i \geq \hat a_{i+1} + b_{i+1}$, this implies that
$m\in {\mathcal P}(\un{a})$.
\QED

\subsection{Combinatorial problems}
The general aim now is to investigate into the structure of 
${\mathcal P}(\un{a})$, for a given sequence $\un{a}\in{\mathbb N}_0^e$:
By (\ref{mainlineprop}) we have  
$$ \mu(\un{a}):=\mathrm {min} ~ {\mathcal P}(\un{a}) 
=\mathrm {min} ~ {\mathcal P}(\un{\tilde a}) =\wp(\un{\tilde a}) ,$$
where $\un{\tilde a}\in{\mathcal N}$ is the associated hull sequence. 
Moreover, by (\ref{cofinprop}) the set 
${\mathcal P}(\un{a})={\mathcal P}(\un{\tilde a})$ 
is a co-finite subset of ${\mathbb N}_0$.
In consequence, the problems associated with $\un{a}$ are to determine
the smallest integer $\sigma(\un{a})\in {\mathbb N}_0$ 
such that all $m\geq\sigma(\un{a})$ are elements of ${\mathcal P}(\un{a})$,
and to determine the gap
$\{\mu(\un{a})+1,\dotsc,\sigma(\un{a})-1\}\setminus {\mathcal P}(\un{a})$.

Note that by the proof of (\ref{cofinprop}) we have 
$\mu(\un{a})\leq\sigma(\un{a})\leq\wp(\un{\hat a})$, where 
$\un{\hat a}$ is the associated $p$-enveloping sequence.
Hence in particular we have shown the following:

\subsection{Theorem}\label{mainlinethm}
Given $\un{a}\in{\mathcal N}$ such that $||\un{a}||\geq p-1$, then we have 
$\mu(\un{a})=\sigma(\un{a})=\wp(\un{a})$, that is the associated 
mainline integers are given as 
${\mathcal P}(\un{a})={\mathbb N}_0+\wp(\un{a})$.
\QED 


\section{\bf Talu's Theorem revisited}\label{talusec}

In this section we develop a conceptual approach to describe the 
smooth epimorphisms onto a given abelian $p$-group.
We first prepare the setting:

\subsection{Abelianisations}\label{fuchsdef}
Let $\Gamma$ be a group with signature, given by the $p$-datum
$(x_1, \dotsc, x_f;h)$, where $h\geq 0$, $f\geq 0$ and $x_f>0$;
note that we are allowing for the case $f=0$, where the $p$-datum 
becomes $(-;h)$. Thus $\Gamma$ is generated by the set
$$ \{a_k,b_k ~:~ 1\leq k\leq h\} \quad \cup \quad
\{c_{ij} ~:~ 1\leq i\leq f,\, 1\leq j\leq x_i\} ,$$
subject to the order relations
$$  c_{ij}^{p^i} = 1, 
\quad\textrm{for}\quad 1 \leq i \leq f 
\quad\textrm{and}\quad 1 \leq j \leq x_i, $$
and the long relation
$$ \prod_{k=1}^{h} [a_k, b_k]\cdot 
   \prod_{i=1}^f\prod_{j=1}^{x_i} c_{ij} = 1 .$$

Let $0\leq f'\leq f$ be defined as follows: 
\[ f':=\left\{ \begin{array}{cl}
0, &\textrm{if }\sum_{i=1}^f x_i\leq 1, \\ 
{\mathrm {max}}\{1\leq d\leq f ~:~ \sum_{i=d}^f x_i\geq 2\}, 
   &\textrm{if }\sum_{i=1}^f x_i\geq 2. \rule{0em}{1.2em} \\
\end{array}\right. \]

In other words, we have $f'=0$ if and only if the $p$-datum is 
$(-;h)$ or $(0,\dotsc,0,1;h)$, 
while otherwise we have $f'=f$ if and only if 
$x_f\geq 2$, and if $x_f=1$ then $1\leq f'<f$ is largest such that $x_{f'}>0$.

It follows from the above presentation
that the abelianisation $H:=\Gamma/[\Gamma,\Gamma]$ of $\Gamma$, 
where $[\Gamma,\Gamma]$ denotes the derived subgroup of $\Gamma$,
can be written as 
$$ H \cong \left\{ 
\begin{array}{ll}
{\mathbb Z}^{2h}, & \textrm{if } f'=0, \\
{\mathbb Z}^{2h} \oplus {\mathbb Z}_{p}^{x_1} \oplus
{\mathbb Z}_{p^2}^{x_2} \oplus \cdots \oplus {\mathbb Z}_{p^f}^{x_f-1},
& \textrm{if } f'=f, \rule{0em}{1.2em} \\
{\mathbb Z}^{2h} \oplus {\mathbb Z}_{p}^{x_1} \oplus
{\mathbb Z}_{p^2}^{x_2} \oplus \cdots \oplus {\mathbb Z}_{p^{f'}}^{x_{f'}},
& \textrm{if } 1\leq f'<f. \rule{0em}{1.2em} \\
\end{array}\right. $$
Indeed, identifying the elements of $\Gamma$ 
with their images under the natural map $\Gamma\longrightarrow H$,
we conclude that $H$ is generated by the set
$$ {\mathcal C} := {\mathcal C}_0 \cup {\mathcal C}_1 
\cup \cdots \cup {\mathcal C}_{f-1} \cup {\mathcal C}_f ,$$
reflecting its decomposition as a direct sum of cyclic subgroups, 
where 
$$ \begin{array}{lcll}
{\mathcal C}_0 &:=& \{a_k,b_k\in H ~:~ 1\leq k\leq h\}, \\
{\mathcal C}_i &:=& \{c_{ij}\in H ~:~ 1\leq j\leq x_i\}, &
\quad \textrm{for} \quad 1\leq i\leq f-1, \rule{0em}{1.2em} \\
{\mathcal C}_f &:=& \{c_{fj}\in H ~:~ 1\leq j\leq x_f-1\}. 
\rule{0em}{1.2em} \\
\end{array} $$

\subsection{Abelian groups}\label{abeliandef}
Let $G$ be a non-trivial abelian $p$-group given by
\[ G \cong {\mathbb Z}_p^{r_1} \oplus {\mathbb Z}_{p^2}^{r_2} 
  \oplus \cdots \oplus {\mathbb Z}_{p^e}^{r_e}, \]
where $e\geq 1$, and $r_i\geq 0$ for $1\leq i\leq e-1$, and $r_e \geq 1$.
Moreover, let 
$$ \{g_{ij} ~:~ 1\leq i\leq e,\, 1\leq j\leq r_i\} $$
be a generating set reflecting the decomposition as a 
direct sum of cyclic subgroups.

Proceeding similarly as above, let $0\leq e'\leq e$
be defined as follows:
\[ e':=\left\{ \begin{array}{cl}
0, &\textrm{if }\sum_{i=1}^e r_i\leq 1, \\ 
{\mathrm {max}}\{1\leq d\leq e ~:~ \sum_{i=d}^e r_i\geq 2\}, 
   &\textrm{if }\sum_{i=1}^e r_i\geq 2. \rule{0em}{1.2em} \\
\end{array}\right. \]
Thus, we have $e'=0$ if and only if $G \cong {\mathbb Z}_{p^e}$
is cyclic, while otherwise we have $e'=e$ if and only if $r_e\geq 2$,
and if $r_e=1$ then $1\leq e'<e$ is largest such that $r_{e'}>0$.

Letting $\Omega_i (G) = \{ g \in G ~:~ g^{p^i} = 1 \}$ 
be the characteristic subgroup of $G$ consisting of all elements 
of order dividing $p^i$, where $0\leq i\leq e$,
we observe that $\Omega_{i-1}(G)$ is a subgroup of index $p$ in $\Omega_i(G)$
if and only if $e'<i\leq e$. 
In other words, using the notions developed in {\cite{ms}},
we have $e'<e$ if and only if $G$ is a group of `GK type', 
in which case $e-e'$ coincides with the length of its `GK series', 
see {\cite[Ex.2.3]{ms}}. 
In view of Kulkarni's Theorem (\ref{kulkarnithm}), and  
the comments in {\cite[Sect.1.1]{ms}}, it is not surprising
that this shows up here in disguised form as well.

For the remainder of this section we keep the notation fixed in
(\ref{fuchsdef}) and (\ref{abeliandef}). Now, 
since any group homomorphism from $\Gamma$ to an abelian group
factors through $H$, from (\ref{fuchs}) we get the following:

\subsection{Proposition}\label{smoothprop}
There is a smooth epimorphism
$\phi : \Gamma \longrightarrow G$ if and only if there is an epimorphism 
$\varphi : H:=\Gamma/[\Gamma,\Gamma] \longrightarrow G$ such that 
$\varphi(c_{ij})$ has order $p^i$, 
for $1 \leq i \leq f$ and $1 \leq j \leq x_i$, and
$\prod_{j=1}^{x_f-1} \varphi(c_{fj})$ has order $p^f$.
\QED

\bigskip

Such an epimorphism $\varphi : H \longrightarrow G$
is also said to be {\bf smooth}. Having this in place, we are prepared
to state a necessary and sufficient arithmetic condition 
when there is a smooth epimorphism $\phi : \Gamma \longrightarrow G$.
By (\ref{smoothprop}) this amounts to give such a condition for a smooth 
epimorphism $\varphi : H \longrightarrow G$, which
is done in (\ref{smoothnecthm}) and (\ref{smoothsuffthm}) 
for necessity and sufficiency, respectively. We call this collection of
statements {\bf Talu's Theorem}, for the following reasons:

We pursue a strategy similar to the one employed in
\cite[La.3.2]{tal} and \cite[Thm.3.3]{tal}, where the statements
of (\ref{smoothnecthm}) and (\ref{smoothsuffthm}) are proven
for the case $p$ odd. Here, we are developing a general approach, 
which covers the case $p=2$ as well, and with which we recover the 
results in \cite{tal} in a more conceptual manner.
In particular, we close a gap in the proof of \cite[Thm.3.3]{tal},
where the element there playing a role similar to the element `$g$' 
in our proof of (\ref{smoothsuffthm}) is incorrectly stated.

\subsection{Theorem}\label{smoothnecthm}
If there exists a smooth epimorphism 
$\varphi : H \longrightarrow G$ then we have $f'=f\leq e$, 
and the following inequalities are fulfilled:
$$ 2h + \sum_{j=i}^f x_j \geq 1 + \sum_{j=i}^e r_j, 
\quad\textrm{for}\quad 1\leq i\leq f,\quad\textrm{and}\quad
   2h \geq \sum_{j=f+1}^e r_j .$$
Moreover, if $p=2$ and $e'<f$, then $x_f$ is even.

\subsection*{Proof}
For $0\leq i\leq e$ let 
$\Omega_i (G) = \{g \in G ~:~ g^{p^i} = 1 \}$ and
$\mho_i (G) = \{ g^{p^i} \in G ~:~ g \in G \}$
be the characteristic subgroups of $G$ consisting of all elements of order 
dividing $p^i$, and of all $p^i$-th powers, respectively. 
In particular $\Omega_1 (G)$ is an ${\mathbb F}_p$-vector space,
where ${\mathbb F}_p$ denotes the field with $p$ elements. 

Now, the existence of the smooth epimorphism $\varphi : H \longrightarrow G$
implies $f'=f\leq e$. We have 
$\mho_e (H)\leq\ker(\varphi)$, thus letting
$$ \tilde H :=H/\mho_e (H)\cong 
{\mathbb Z}_{p}^{x_1} \oplus
{\mathbb Z}_{p^2}^{x_2} \oplus \cdots \oplus {\mathbb Z}_{p^f}^{x_f-1} 
\oplus {\mathbb Z}_{p^e}^{2h} $$
yields an epimorphism $\tilde\varphi : \tilde{H} \longrightarrow G$.
Hence dualising we get a monomorphism
$\tilde\varphi^\ast : G^\ast:=\mathrm{Hom}(G,\C^\ast)
\longrightarrow \mathrm{Hom}(\tilde H,\C^\ast)=\tilde H^\ast$,
that is $G\cong G^\ast$ is isomorphic to a subgroup of 
$\tilde H^\ast\cong \tilde H$.
Thus $\Omega_i (G)$ and $\mho_i (G)$ can be identified with subgroups
of $\Omega_i (\tilde H)$ and $\mho_i (\tilde H)$, respectively, and
hence we have
$$\dim_{{\mathbb F}_p}(\Omega_1 (\mho_i (G))) \leq
 \dim_{{\mathbb F}_p}(\Omega_1 (\mho_i (\tilde H))) .$$
  
Now, for $0\leq i\leq e-1$ we have
$$ \Omega_1 (\mho_i (G))\cong {\mathbb Z}_p^{r_{i+1}} \oplus 
 {\mathbb Z}_p^{r_{i+2}} \oplus \cdots \oplus {\mathbb Z}_p^{r_e}, $$
which yields
$$ \dim_{{\mathbb F}_p}(\Omega_1 (\mho_i (G)))=\sum_{j=i+1}^e r_j .$$
Similarly, for $0\leq i\leq f-1$ we have
$$ \Omega_1 (\mho_i (\tilde H))\cong {\mathbb Z}_p^{x_{i+1}} \oplus 
 {\mathbb Z}_p^{x_{i+2}} \oplus \cdots \oplus {\mathbb Z}_p^{x_{f-1}}
 \oplus {\mathbb Z}_p^{x_f-1} \oplus {\mathbb Z}_p^{2h} ,$$
yielding
$$ \dim_{{\mathbb F}_p}(\Omega_1 (\mho_i ( \tilde H)))=
2h-1+\sum_{j=i+1}^f x_j ,$$
while for $f\leq i\leq e-1$  we get
$$ \dim_{{\mathbb F}_p}(\Omega_1 (\mho_i ( \tilde H)))=2h .$$
 
Finally, let $p=2$ and $e'<f\leq e$. Then $G$ has shape
$$ G \cong {\mathbb Z}_{2}^{r_1} \oplus {\mathbb Z}_{4}^{r_2} \oplus 
  \dotsc \oplus {\mathbb Z}_{2^{e'}}^{r_{e'}} \oplus {\mathbb Z}_{2^e} ,$$
and thus 
$$ \Omega_f(G)/\Omega_{f-1}(G)\cong 
{\mathbb Z}_{2^{e-f}}/ {\mathbb Z}_{2^{e-f+1}} \cong {\mathbb Z}_2 .$$
Now we observe that 
$\varphi(c_{fj})\in\Omega_f(G)\setminus\Omega_{f-1}(G)$,
for $1\leq j\leq x_f-1$, where 
$\prod_{j=1}^{x_f-1}\varphi(c_{fj})\not\in\Omega_{f-1}(G)$
as well, implying that $x_f-1$ is odd.
\QED

\subsection{Theorem}\label{smoothsuffthm}
Let $f'=f\leq e$, where in case $p=2$ and $e'<f$ we
additionally assume that $x_f$ is even, such that 
$$ 2h + \sum_{j=i}^f x_j \geq 1 + \sum_{j=i}^e r_j,
\quad\textrm{for}\quad 1\leq i\leq f,\quad\textrm{and}\quad
   2h \geq \sum_{j=f+1}^e r_j .$$
Then there exists a smooth epimorphism 
$\varphi : H \longrightarrow G$.

\subsection*{Proof}
By the inequalities assumed we have
$$ |{\mathcal C}_0\cup{\mathcal C}_f\cup{\mathcal C}_{f-1}
\cup\cdots\cup{\mathcal C}_i| \geq \sum_{j=i}^e r_j,
\quad\textrm{for}\quad 1\leq i\leq f ,
\quad\textrm{and}\quad |{\mathcal C}_0| \geq \sum_{j=f+1}^e r_j ,$$
where the latter sum is empty if $e=f$. 
Thus we may choose a subset
${\mathcal D}_{f+1}\subseteq {\mathcal C}_0$
of cardinality $\sum_{j=f+1}^e r_j$.
Subsequently, for $f\geq i\geq 1$ we may recursively choose,
disjointly from ${\mathcal D}_{f+1}$, pairwise disjoint sets 
$$ {\mathcal D}_i=\{d_{i,1},\dotsc,d_{i,r_i}\}\subseteq 
{\mathcal C}_0\cup{\mathcal C}_f\cup{\mathcal C}_{f-1}
\cup\cdots\cup{\mathcal C}_i $$
of cardinality $r_i$. Let
$$ {\mathcal C}'_i:={\mathcal C}_i\setminus
   \bigg(\bigcup_{j=1}^i{\mathcal D}_j\bigg)
\quad\textrm{for}\quad 1\leq i\leq f, \quad\textrm{and}\quad 
   {\mathcal C}'_0:={\mathcal C}_0\setminus
   \bigg(\bigcup_{j=1}^{f+1}{\mathcal D}_j\bigg) .$$
We are going define a homomorphism $\varphi : H\longrightarrow G$
by specifying the image of ${\mathcal C}$:

The direct summand $\langle{\mathcal D}_{f+1}\rangle$
of $H$ is a free abelian group of rank $\sum_{j=f+1}^e r_j$, 
hence choosing $\varphi(c)$ appropriately, 
for $c\in {\mathcal D}_{f+1}\subseteq{\mathcal C}_0$,
the direct summand
\[ G':=\langle g_{ij} ~:~ f+1\leq i\leq e,\, 1 \leq j\leq r_i \rangle
\cong {\mathbb Z}_{p^{f+1}}^{r_{f+1}} \oplus 
   {\mathbb Z}_{p^{f+2}}^{r_{f+2}} \oplus \cdots \oplus 
{\mathbb Z}_{p^e}^{r_e} \]
of $G$ becomes an epimorphic image of $\langle{\mathcal D}_{f+1}\rangle$.
Thus letting $\varphi(c):=1$ for $c\in{\mathcal C}'_0$,
we are done in the case $f=0$. Hence we may assume that 
$f'=f>0$, thus we have $x_f\geq 2$ and ${\mathcal C}_f\neq\emptyset$,
where we may assume that
${\mathcal C}_f \cap{\mathcal D}_f\neq\emptyset$ whenever $r_f>0$.

Now, for $d_{ij}\in{\mathcal C}_0\cap {\mathcal D}_i$, 
where $1\leq i\leq f$, we let $\varphi(d_{ij}):=g_{ij}$.
Moreover, for $d_{ij}\in{\mathcal C}_k\cap{\mathcal D}_i$,
where $1\leq i\leq k<f\leq e$, we let
$\varphi(d_{ij}):=g_{ij}\cdot g_{e,r_e}^{p^{e-k}}$,
while for $c\in{\mathcal C}'_k$ we let $\varphi(c):=g_{e,r_e}^{p^{e-k}}$.
To specify $\varphi(c)$ for $c\in {\mathcal C}_f$ we need some flexibility: 

For $d_{ij}\in{\mathcal C}_f\cap{\mathcal D}_i$,
where $1\leq i\leq f$, we let $\varphi(d_{ij})=g_{ij}\cdot c'$, 
for some $c'\in G$, while for $c\in{\mathcal C}'_f$ we just
write $\varphi(c)=c'$.
Then we have to show that the elements $c'$ can be chosen suitably
to give rise to an epimorphism such that all $\varphi(c)$, 
where $c\in{\mathcal C}_f$, as well as
$g:=\prod_{c\in{\mathcal C}_f}\varphi(c)$ have order $p^f$.

In particular, $\varphi(c)$ will have order $p^f$, if 
$c\in{\mathcal C}_f\setminus{\mathcal D}_f$ and 
$c'\in G$ is chosen to have order $p^f$,
or if $c\in{\mathcal C}_f\cap{\mathcal D}_f$ and 
$c'\in G'$ is chosen to have order dividing $p^f$.
Moreover, $\varphi$ will be an epimorphism whenever $f<e$ and we choose
$c'\in G'$ for all $c\in{\mathcal C}_f\cap(\bigcup_{i=1}^f{\mathcal D}_i)$.
The order condition on $g$ will be checked by showing 
that the image of $g$ under a suitable projection of $G$ onto one
of its direct summands already has order $p^f$.
We now distinguish various cases:

{\bf i)}
Let $f<e'\leq e$. Then pick $c_0\in{\mathcal C}_f$, and let
$c'_0:=g_{e',1}^{p^{e'-f}}$, while for $c_0\neq c\in{\mathcal C}_f$ 
let $c':=g_{e,r_e}^{p^{e-f}}$;
note that for $e'=e$ we have $r_e\geq 2$. 
Then projecting $g$ onto $\langle g_{e',1}\rangle$ yields $c'_0$,
which has order $p^f$.

{\bf ii)}
Let $f=e'\leq e$. Then, since $r_f=r_{e'}>0$, we may assume that
$d_{e',1}\in{\mathcal C}_f\cap {\mathcal D}_f$.
For $c\in{\mathcal C}_f\setminus{\mathcal D}_f$ let 
$c':=g_{e,r_e}^{p^{e-f}}$, while for 
$c\in {\mathcal C}_f\cap {\mathcal D}_f$ let $c':=1$;
note that for $f=e'=e$ we have $r_e\geq 2$,
and $d_{e,r_e}\in {\mathcal C}_0\cup {\mathcal C}_f$ 
implies that $\varphi$ is an epimorphism.
Projecting $g$ onto $\langle g_{e',1}\rangle$ yields $g_{e',1}$,
which has order $p^f$.

{\bf iii)}
Let $e'<f<e$. Then for $c\in{\mathcal C}_f$ let 
$c':=(g_{e,1}^{p^{e-f}})^{a_c}$, where $a_c$ is chosen coprime to $p$.
Projecting $g$ onto $\langle g_{e,1}\rangle$ yields 
$(g_{e,1}^{p^{e-f}})^a$, where $a:=\sum_{c\in{\mathcal C}_f}a_c$.
The latter element has order $p^f$ if and only if $a$ is coprime to $p$.
If $p$ is odd, this can be achieved by picking any $c\in{\mathcal C}_f$
and replacing $a_c$ by $a_c+1$ or $a_c-1$, if necessary.
If $p=2$, then $a_c$ is odd for all $c\in{\mathcal C}_f$, which, 
since $|{\mathcal C}_f|=x_f-1$ is odd,
implies that $a$ is odd. 

{\bf iv)}
Let $e'<f=e$.
Then, since $r_f=r_e=1$, we may assume that
${\mathcal C}_f\cap {\mathcal D}_f=\{d_{e,1}\}$.
For $c\in{\mathcal C}_f$ let $c':=g_{e,1}^{a_c}$,
where $a_c$ is chosen coprime to $p$ for $c\neq d_{e,1}$, while 
for $c=d_{e,1}$ we choose $a_c$ such that $1+a_c$ is coprime to $p$.
This implies that $\varphi(d_{e,1})$ has order $p^f$ and that $\varphi$ 
is an epimorphism.
Projecting $g$ onto $\langle g_{e,1}\rangle$ yields 
$g_{e,1}^a$, where $a:=1+\sum_{c\in{\mathcal C}_f}a_c$.
The latter element has order $p^f$ if and only if $a$ is coprime to $p$.
If $p$ is odd, this can be achieved by picking $c\in{\mathcal C}_f$
and replacing $a_c$ by $a_c+1$ or $a_c-1$, if necessary.
If $p=2$, then $a_c$ is odd for all $d_{e,1}\neq c\in{\mathcal C}_f$, 
and $1+a_c$ is odd for $c=d_{e,1}$, which, 
since $|{\mathcal C}_f|=x_f-1$ is odd,
implies that $a$ is odd. 
\QED

\section{\bf Transforming to mainline integers}\label{transsec}

In this section we show how mainline integers, as introduced in
Section \ref{mainlinesec}, can be reconciled with the problem of
determining the (reduced) genus spectrum of abelian $p$-groups 
and the results of Section \ref{talusec}.

\subsection{Translating the reduced genus map}\label{preprem} 
Let still $G$ be a non-trivial abelian $p$-group of exponent $p^e$.

{\bf i)}
In order to reformulate the results of Section \ref{talusec},
we define 
$\alpha : {\mathbb N}_0^{e+1} \longrightarrow {\mathbb N}_0^{e+1}$ by 
\[ \alpha(x_1, \dotsc, x_e;x_0) := 
\bigg(\sum_{i=1}^{e} x_i+2x_0, \sum_{i=2}^{e} x_i+2x_0, \dotsc,
      x_e+2x_0, 2x_0\bigg), \]
which is injective and has image, 
using the notation from (\ref{mainlinesetting}), 
\[ {\mathrm {im}}(\alpha) = {\mathcal N}'(e+1) :=
\{ (a_1, \dotsc, a_{e+1}) \in {\mathcal N}(e+1) ~:~ a_{e+1} \in 2\N_0 \}. \]
The inverse map 
$\alpha^{-1} : {\mathcal N}'(e+1) \longrightarrow {\mathbb N}_0^{e+1}$ 
is given by 
\[ \alpha^{-1} (a_1, \dotsc, a_{e+1}) 
:= (a_1 - a_2, \dotsc, a_e - a_{e+1}; \frac{a_{e+1}}{2}) .\]

Letting $D(G) \subset {\mathbb N}_0^{e+1}$ be the data spectrum of $G$ 
as introduced in (\ref{genusmap}), let
$$ A(G):=\alpha(D(G)) \subset {\mathbb N}_0^{e+1} .$$ 
Then the reduced genus map
$g_0 : D(G) \longrightarrow \frac{1}{2}\cdot(\{-1\}\cup\N_0)$,
given by
\[ g_0(x_1, \dotsc, x_e;h) = 
-p^e + \bigg(h+\frac {1}{2}\cdot\sum_{i=1}^{e} x_i\bigg)\cdot p^e
- {\frac {1}{2}} \cdot\sum_{i=1}^{e} x_ip^{e-i} ,\]
can be rephrased as
$\gamma = g_0 \circ \alpha^{-1} : A(G) 
\longrightarrow \frac{1}{2}\cdot(\{-1\}\cup\N_0)$, where explicitly
\[ \gamma(a_1, \dotsc, a_{e+1}) = - p^e + 
{\frac {a_{e+1}}{2}} + {\frac {p-1}{2}}\cdot\wp(a_1, \dotsc, a_e) .\]

{\bf ii)}
As will become clear below, elements of the form
$(x_1, \dotsc, x_i,0,\dotsc,0;h)\in D(G)$, for some  $0 \leq i \leq e$,
are of particular importance. These translate into elements of the form
$(a_1, \dotsc, a_i, 2a, \dotsc, 2a) \in {\mathcal N}'(e+1)$. For the
latter we have 
$$ \gamma(a_1, \dotsc, a_i, 2a, \dotsc, 2a)
= -p^e + a + {\frac {p-1}{2}} \cdot \wp(a_1, \dotsc, a_i, 2a, \dotsc, 2a) ,$$
where the argument of $\wp$ is a sequence of length $e$, and yields
$$\wp(a_1, \dotsc, a_i, 2a, \dotsc, 2a) =
p^{e-i}\cdot\sum_{j=1}^i a_jp^{i-j} + 2a \cdot \sum_{j=0}^{e-i-1} p^j .$$
From that we get
$$ \gamma(a_1, \dotsc, a_i, 2a, \dotsc, 2a)
= -p^e + p^{e-i}\cdot\bigg(a+\frac{p-1}{2}\cdot\wp(a_1,\dotsc,a_i)\bigg). $$

In particular, for $i=0$ we get
$\gamma(2a, \dotsc, 2a)= (a-1)\cdot p^e$,
while for $i=e$ we recover 
$\gamma(a_1, \dotsc, a_e,2a)=-p^e + a +
{\frac {p-1}{2}} \cdot \wp(a_1, \dotsc, a_e)$.
Note that we have $\gamma(a_1, \dotsc, a_i, 2a, \dotsc, 2a)\in\Z$,
unless $p=2$ and $i=e$ and $a_e$ odd, in which case we have
$\gamma(a_1, \dotsc, a_e, 2a)\in\frac{1}{2}\Z\setminus\Z$.

\subsection{Translating Talu's Theorem}\label{talurem} 
Let again
$G \cong {\mathbb Z}_p^{r_1} \oplus {\mathbb Z}_{p^2}^{r_2} 
  \oplus \cdots \oplus {\mathbb Z}_{p^e}^{r_e} $,
where $e\geq 1$, and $r_i\geq 0$ for $1\leq i\leq e-1$, and $r_e \geq 1$. 
Moreover, for $1 \leq i \leq e+1$ we fix 
\[ s_i := 1+\sum_{j=i}^e r_j, \] 
Hence we have $\un{s}:=(s_1,\dotsc,s_{e+1})\in{\mathcal N}(e+1)$ 
such that $s_e\geq 2$ and $s_{e+1}=1$.
Having this in place, (\ref{smoothnecthm}) and (\ref{smoothsuffthm}) 
can be rephrased as follows: 

{\bf i)}
For $p$ odd we have
\[ A(G) := A_0 \cup A_1 \cup \cdots \cup A_e , \]
where for $0 \leq i \leq e$ we let, setting $a_0:=\infty$, 
$$ \begin{array}{rcl} 
A_i & := &\{ \un{a}\in {\mathcal N}'(e+1) ~:~ 
(a_1,\dotsc,a_i) \geq (s_1, \dotsc, s_i), \\
&& \quad 
a_{i+1} = \cdots = a_{e+1} \geq s_{i+1} - 1,\, a_i - a_{i+1} \geq 2 \}.
\rule{0em}{1.2em} \\
\end{array} $$
In particular, we have
$$ A_0 = \{ \un{a}\in {\mathcal N}'(e+1) ~:~ 
a_1 = \cdots = a_{e+1} \geq s_1 - 1 \} $$
and
\[ A_e = \{ \un{a}\in {\mathcal N}'(e+1) ~:~
(a_1,\dotsc,a_e) \geq (s_1, \dotsc, s_e),\, a_e - a_{e+1} \geq 2 \}. \]

For $0\leq i<j\leq e$ the sequences in $A_i$ satisfy $a_j = a_{j+1}$,
while those in $A_j$ satisfy $a_j - a_{j+1} \geq 2$, 
hence $A_i \cap A_j = \emptyset$, thus
$A(G)$ is disjointly covered by the $A_i$.

{\bf ii)}
For $p=2$, letting $0\leq e'\leq e$ be as defined in (\ref{abeliandef}), 
we get
\[ A(G) := A_0 \cup A_1 \cup \cdots \cup A_{e'}
\cup A'_{e'+1} \cup \cdots \cup A'_e , \]
where for $1 \leq i \leq e$ we let
$$ A'_i := \{\un{a}\in A_i ~:~ a_i - a_{i+1} \in 2\N \} .$$
In particular, for $i=e$ we get
$$ A'_e := \{\un{a}\in A_e ~:~ a_e \in 2\N \} .$$
Note that we have $\gamma(A_e)\subseteq\frac{1}{2}\Z$ 
and $\gamma(A'_e)\subseteq\Z$, thus we recover Kulkarni's Theorem
(\ref{kulkarnithm}) in the case of abelian $p$-groups.

\subsection{Towards the minimum genus}\label{minrem}
This now gives a handle to compute the reduced minimum genus of $G$, 
which for $p$ odd is given as
$$ \mu_0(G) = 
\mathrm {min}\{ {\mathrm {min}} ~\gamma (A_i) ~:~ 0 \leq i\leq e \} ,$$
while for $p=2$ we get
$$ \mu_0(G) = 
\mathrm {min}~( \{ {\mathrm {min}} ~\gamma (A_i) ~:~ 0 \leq i\leq e' \}
\cup \{ {\mathrm {min}} ~\gamma (A'_i) ~:~ e'< i\leq e \} ) .$$

{\bf i)}
We proceed to derive formulae, in terms of the sequence 
$\un{s}=(s_1,\dotsc,s_{e+1})$ associated with $G$, 
to determine ${\mathrm {min}} ~\gamma (A_i)$, 
for $0\leq i\leq e$: To this end, let
$$\un{s}^i:= (s_1,\dotsc,s_i,
2\cdot \lfloor {\frac {s_{i+1}}{2}} \rfloor, \dotsc, 
2\cdot \lfloor {\frac {s_{i+1}}{2}} \rfloor) \in {\mathcal N}'(e+1) $$
and
$$ \un{s}^{i+}:=(s_1,\dotsc,s_{i-1},s_i+\epsilon_i,
2\cdot\lfloor {\frac {s_{i+1}}{2}} \rfloor,\dotsc,
2\cdot\lfloor {\frac {s_{i+1}}{2}} \rfloor) \in {\mathcal N}'(e+1) ,$$
where $\epsilon_i\in\{0,1,2\}$ is chosen minimal such that
$s_i+\epsilon_i-2\cdot \lfloor {\frac {s_{i+1}}{2}} \rfloor \geq 2$,
that is
$$ \epsilon_i:=\left\{ 
\begin{array}{rl}
0, & \textrm{if } s_i-s_{i+1}\geq 2, \\
0, & \textrm{if } s_i-s_{i+1}=1,\, s_{i+1} \textrm{ odd}, \rule{0em}{1.2em} \\
1, & \textrm{if } s_i-s_{i+1}=1,\, s_{i+1} \textrm{ even}, \rule{0em}{1.2em} \\
1, & \textrm{if } s_i=s_{i+1},\, s_{i+1} \textrm{ odd}, \rule{0em}{1.2em} \\
2, & \textrm{if } s_i=s_{i+1},\, s_{i+1} \textrm{ even}. \rule{0em}{1.2em} \\
\end{array} \right. $$
Note that for $i=e$ we have $s_{e+1}=1$ and $s_e\geq 2$, 
and thus $\epsilon_e=0$; moreover, for $i=0$ we let $\epsilon_0=0$.

It now follows from the description of $A_i$, and (\ref{mainlineprop}), that 
${\mathrm {min}} ~\gamma (A_i)$ is attained precisely for the hull sequence
$$ \tilde{\un{s}}^{i+}=(\tilde s_1,\dotsc,\tilde s_i,
2\cdot \lfloor {\frac {s_{i+1}}{2}} \rfloor, \dotsc,
2\cdot \lfloor {\frac {s_{i+1}}{2}} \rfloor) \in {\mathcal N}'(e+1) ,$$
of $\un{s}^{i+}$, where 
the prefix $(\tilde s_1,\dotsc,\tilde s_i)$ of length $i$ is
determined as follows:

For $i\geq 1$ let $0\leq i''\leq i'<i$ be both maximal such that 
$s_{i'}-s_i\geq 1$ and $s_{i''}-s_i\geq 2$;
hence, if $i''<i'$ then we have $s_{i'}-s_{i'+1}=1$, and 
$i'=0$ and $i''=0$ refer to the cases $s_1=s_i$ and $s_1-s_i\leq 1$,
respectively. 
Then $(\tilde s_1,\dotsc,\tilde s_i)$ is given as
$$ \begin{array}{cl}
(s_1,\dotsc,s_i), & \textrm{if } \epsilon_i=0, \\
(s_1,\dotsc,s_{i'},s_{i'+1}+1,\dotsc,s_i+1), & \textrm{if } \epsilon_i=1, 
\rule{0em}{1.2em} \\
(s_1,\dotsc,s_{i''},s_{i''+1}+1,\dotsc,s_{i'}+1,s_{i'+1}+2,\dotsc,s_i+2),
& \textrm{if } \epsilon_i=2. \rule{0em}{1.2em} \\
\end{array} $$

Thus letting 
$$ \mu_i := \gamma(\un{s}^i) =
-p^e + p^{e-i}\cdot \bigg(\lfloor {\frac {s_{i+1}}{2}} \rfloor + 
{\frac {p-1}{2}} \cdot \wp(s_1,\dotsc,s_i) \bigg), $$
we get 
$$ {\mathrm {min}} ~\gamma (A_i)
=\gamma (\tilde{\un{s}}^{i+})
=\left\{ \begin{array}{ll}
\mu_i, & \textrm{if } \epsilon_i=0,  \\
\mu_i+\frac{1}{2} \cdot (p^{e-i'}-p^{e-i}), & \textrm{if } \epsilon_i=1, 
\rule{0em}{1.2em} \\
\mu_i+\frac{1}{2} \cdot (p^{e-i''}+p^{e-i'}-2p^{e-i}), & 
\textrm{if } \epsilon_i=2. \rule{0em}{1.2em} \\
\end{array} \right.  $$

In particular, we have 
$$ {\mathrm {min}} ~ \gamma (A_e) = \mu_e = 
-p^e + \frac {p-1}{2} \cdot \wp(s_1,\dotsc,s_e), $$
being attained precisely for $(s_1, \dotsc, s_e, 0)$,
and 
$$ {\mathrm {min}} ~ \gamma (A_0) = \mu_0 
= (\lfloor {\frac {s_1}{2}} \rfloor-1)\cdot p^e ,$$
being attained precisely for 
$(2\cdot \lfloor {\frac {s_1}{2}} \rfloor, \dotsc, 
2\cdot \lfloor {\frac {s_1}{2}} \rfloor)$.

{\bf ii)}
It remains to consider ${\mathrm {min}} ~\gamma (A'_i)$, for $e'<i\leq e$,
in the case $p=2$: For $e'<i<e$ we have $s_i=s_{i+1}=2$, hence 
$\tilde s_i=4$ and $2\cdot\lfloor {\frac {s_{i+1}}{2}} \rfloor =2$,
while for $e'<i=e$ we have $s_e=2$ and $s_{e+1}=1$, hence 
$\tilde s_e=2$ and $2\cdot\lfloor {\frac {s_{e+1}}{2}} \rfloor =0$.
Thus the above description for $e'<i\leq e$ yields
$$ {\mathrm {min}} ~\gamma (A'_i) = {\mathrm {min}} ~\gamma (A_i) 
= \gamma (\tilde{\un{s}}^{i+}) ,$$
implying that the reduced minimum genus of $G$, just as for $p$ odd,
is given as 
$$ \mu_0(G) = 
\mathrm {min}\{ {\mathrm {min}} ~\gamma (A_i) ~:~ 0 \leq i\leq e \} .$$

\subsection{Further towards the minimum genus}\label{minbetterrem}
We turn to the question whether there are relations between the various 
$\gamma (\tilde{\un{s}}^{i+})={\mathrm {min}} ~\gamma (A_i)$, 
for $0\leq i\leq e$, which would allow to take the minimum
determining $\mu_0(G)$ over a smaller set. To this end, we consider
the cases where $\epsilon_i\neq 0$; hence we have $1\leq i\leq e-1$:  

{\bf i)}
If $s_{i+1}$ is even and $s_i=s_{i+1}$, then we have
$$ \begin{array}{lcl}
\un{s}^{i+}&=&(s_1,\dotsc,s_{i-1},s_i+2,s_i,\dotsc,s_i), \\
\un{s}^{(i-1)+}&=&(s_1,\dotsc,s_{i-1}+\epsilon_{i-1},s_i,s_i,\dotsc,s_i), 
\rule{0em}{1.2em} \\
\end{array} $$ 
where $\epsilon_{i-1}=0$ whenever $s_{i-1}\geq s_i+2$, and
$s_{i-1}+\epsilon_{i-1}=s_i+2$ otherwise.

{\bf ii)}
If $s_{i+1}$ is even and $s_i-s_{i+1}=1$, then we have
$$ \begin{array}{lcl}
\un{s}^{i+}&=&(s_1,\dotsc,s_{i-1},s_i+1,s_i-1,\dotsc,s_i-1), \\ 
\un{s}^{(i-1)+}&=&
(s_1,\dotsc,s_{i-1}+\epsilon_{i-1},s_i-1,s_i-1,\dotsc,s_i-1), 
\rule{0em}{1.2em} \\
\end{array} $$ 
where $\epsilon_{i-1}=0$ whenever $s_{i-1}\geq s_i+1$, and
$s_{i-1}+\epsilon_{i-1}=s_i+1$ otherwise.

{\bf iii)}
If $s_{i+1}$ is odd and $s_i=s_{i+1}$, then we have
$$ \begin{array}{lcl}
\un{s}^{i+}&=&(s_1,\dotsc,s_{i-1},s_i+1,s_i-1,\dotsc,s_i-1), \\
\un{s}^{(i-1)+}&=&
(s_1,\dotsc,s_{i-1}+\epsilon_{i-1},s_i-1,s_i-1,\dotsc,s_i-1), 
\rule{0em}{1.2em} \\
\end{array} $$ 
where $\epsilon_{i-1}=0$ whenever $s_{i-1}\geq s_i+1$, and
$s_{i-1}+\epsilon_{i-1}=s_i+1$ otherwise.
\QED

\bigskip 

Hence, in either of these cases, going over to hull sequences yields
$\tilde{\un{s}}^{i+} \geq \tilde{\un{s}}^{(i-1)+}$, implying
${\mathrm {min}} ~\gamma (A_i)=\gamma(\tilde{\un{s}}^{i+} )\geq
\gamma(\tilde{\un{s}}^{(i-1)+})={\mathrm {min}} ~\gamma (A_{i-1})$.
Thus ${\mathrm {min}} ~\gamma (A_i)$ need not be considered in 
finding $\mu_0(G)$. 
Hence we are left with the cases $0\leq i\leq e$
such that $\epsilon_i=0$, that is ${\mathrm {min}} ~\gamma (A_i)=\mu_i$.

Moreover, if $s_1$ is even, then since $s_1\geq\cdots\geq s_e\geq 2$ we have
$$ \mu_e = -p^e + \frac {p-1}{2} \cdot \wp(s_1,\dotsc,s_e)
\leq -p^e + \frac {s_1}{2} \cdot (p^e-1)
<(\frac{s_1}{2}-1)\cdot p^e =\mu_0  ,$$
hence in this case ${\mathrm {min}} ~\gamma (A_0)$ need not be considered
in finding $\mu_0(G)$.
Thus, in conclusion, we have proved the following:

\subsection{Theorem}\label{minthm}
Keeping the above notation, we have
$$ \mu_0(G) = 
{\mathrm {min}}\{{\mathrm {min}} ~\gamma (A_i) ~:~ i\in \mathcal I(G) \}
= {\mathrm {min}}\{ \mu_i~:~ i\in \mathcal I(G) \}, $$
where, letting $s_0:=\infty$, we have
$$ \mathcal I(G):= 
\{0\leq i\leq e~:~ s_i-s_{i+1}\geq 2\} \quad \cup \quad 
\{0\leq i\leq e~:~ s_i-s_{i+1}=1,\,s_{i+1}\textrm{ odd} \} .$$
In particular, we always have $\{0,e\}\subseteq \mathcal I(G)$, 
but if $s_1$ is even then to find $\mu_0(G)$ it suffices 
to consider $i\in\mathcal I(G)\setminus\{0\}$ only
\QED

\bigskip

In other words, finding $\mu_0(G)$ is reduced to computing the minimum of
$|\mathcal I(G)|\leq e+1$ numbers, which are given explicitly in terms of 
known invariants of $G$. 
In particular, this machinery to determine $\mu_0(G)$ will feature 
prominently in the proof of our main result (\ref{genusthm}).
Moreover, to underline the effectiveness of these techniques,
in Sections \ref{smallranksec} and \ref{smallexpsec} we give 
detailed example treatments of the groups of rank at most $2$,
and of the groups of exponent at most $p^2$, respectively.

\subsection{Translating back}\label{backrem}
We translate the results back, 
to express $\mu_i={\mathrm {min}} ~\gamma (A_i)$,
for $i\in \mathcal I(G)$, in terms of the $p$-datum giving rise to $\mu_i$,
which by (\ref{preprem}) is given as
$$ \un{x}^i = (x_1,\dotsc,x_e;h):= 
\alpha^{-1}(s_1,\dotsc,s_i,
2\cdot \lfloor {\frac {s_{i+1}}{2}} \rfloor, \dotsc, 
2\cdot \lfloor {\frac {s_{i+1}}{2}} \rfloor) .$$

{\bf i)}
If $r_i=s_i-s_{i+1}\geq 2$ and $s_{i+1}$ is even, then we have 
$$ \un{x}^i =
(r_1,\dotsc,r_i,0,\dotsc,0; \frac {s_{i+1}}{2} ) ,$$
yielding
$$ \mu_i=
p^e\cdot\bigg( {\frac {s_{i+1}}{2}} -1+
\frac{1}{2}\cdot\sum_{j=1}^i r_j(1-\frac{1}{p^j}) \bigg) .$$

{\bf ii)}
If $r_i=s_i-s_{i+1}\geq 1$ and $s_{i+1}$ is odd, then we have 
$$ \un{x}^i =
(r_1,\dotsc,r_{i-1},r_i+1,0,\dotsc,0;\frac {s_{i+1}-1}{2}) ,$$
yielding
$$ \mu_i=
p^e\cdot\bigg(\frac {s_{i+1}-1}{2}-1+
\frac{1}{2}\cdot\sum_{j=1}^i r_j(1-\frac{1}{p^j})+
\frac{1}{2}\cdot(1-\frac{1}{p^i})\bigg) .$$

In particular, the case $i=0$ is encompassed by the above cases, 
depending on whether $s_1$ is even or odd, respectively, by
$\un{x}^0 = (0,\dotsc,0; \lfloor {\frac {s_1}{2}} \rfloor)$,
where this case need not be considered if $s_1$ is even.
Moreover, the case $i=e$, since $s_{e+1}=1$, is subsumed in the second 
of the above cases, by $\un{x}^e = (r_1,\dotsc,r_{e-1},r_e+1;0)$.

Finally, the various $\mu_i={\mathrm {min}} ~\gamma(A_i)$ to 
be considered belong to pairwise distinct orbit genera, inasmuch the map 
$$ \mathcal I(G)\longrightarrow \Z: 
   i\mapsto \lfloor {\frac {s_{i+1}}{2}} \rfloor $$
is strictly decreasing, hence in particular is injective: 
Indeed, if $i-1,i\in \mathcal I(G)$, then we have $s_i-s_{i+1}\geq 1$ anyway;
and if $s_i$ is odd and $s_{i+1}$ is even, then from $s_i-s_{i+1}\geq 2$ 
we still get  
$\lfloor {\frac {s_i}{2}} \rfloor =\frac {s_i-1}{2}>\frac {s_{i+1}}{2}
=\lfloor {\frac {s_{i+1}}{2}} \rfloor$. 

\subsection{Maclachlan's method}\label{macrem}
We compare our approach with the method to compute the minimum
genus for arbitrary non-cyclic abelian groups given in \cite{mac}:

Let $G$ be a non-cyclic abelian group, with sequence 
$(n_1,\dotsc,n_s)$ of invariants giving rise to
the Smith normal form abelian group presentation of $G$;
hence we have $s\geq 2$, and the exponent of $G$ equals $n_s$.
Let $\nu_h\in\N_0$ be the reduced minimum genus 
afforded by all signatures of $G$ with fixed orbit genus $h\geq 0$.
Then, by \cite[Thm.4]{mac}, the reduced minimum genus of $G$ equals
$$ \mu_0(G)=
{\mathrm {min}} \{\nu_h ~:~ 0\leq h\leq \lfloor\frac{s}{2}\rfloor\} ,$$
where the numbers $\nu_h$ can be computed explicitly as
$$ \nu_h=n_s\cdot\bigg(h-1+
\frac{1}{2}\cdot\sum_{k=1}^{s-2h}(1-\frac{1}{n_k})+ 
\frac{1}{2}\cdot(1-\frac{1}{n_{s-2h}}) \bigg) .$$

In our case of abelian $p$-groups this reads as follows: We have
$$ (n_1,\dotsc,n_s)=
   (p,\dotsc,p,\,p^2,\dotsc,p^2,\,\dotsc,\,p^e,\dotsc,p^e) ,$$ 
where the entry $p^i$ occurs $r_i$ times, for $1\leq i\leq e$;
hence we have $s=\sum_{i=1}^e r_i=s_1-1$.
Thus we are able to improve \cite[Thm.4]{mac}, 
for non-cyclic abelian $p$-groups, as follows:
By the injectivity of the map
$\mathcal I(G)\longrightarrow \Z: 
 i\mapsto \lfloor {\frac {s_{i+1}}{2}} \rfloor$, 
for $i\in\mathcal I(G)$ we have 
$$ \nu_{\lfloor {\frac {s_{i+1}}{2}} \rfloor} = \mu_i ,$$
and thus by (\ref{minthm}) we may compute $\mu_0(G)$
as a minimum over a set of cardinality $|\mathcal I(G)|\leq e+1$ 
instead of one of cardinality $\lfloor\frac{s_1-1}{2}\rfloor+1$, as
$$ \mu_0(G)= {\mathrm {min}} 
\{\nu_{\lfloor {\frac {s_{i+1}}{2}} \rfloor} ~:~ i\in\mathcal I(G)\} .$$

Recall that whenever $s_1$ is even the case $i=0$ need not be considered,
so that we always get a subset of the indices used in \cite{mac}.
From the formulae in (\ref{backrem}) to compute $\mu_i$
in terms of $p$-data, we recover the formulae for 
$\nu_{\lfloor {\frac {s_{i+1}}{2}} \rfloor}$ given in \cite{mac}.
Finally, we point out that our approach is also valid for cyclic
$p$-groups, while cyclic groups are excluded in \cite{mac}.
Moreover, since in \cite{mac} only genera $g\geq 2$ are
considered, the case $s=2$ and some small abelian groups
have to be treated as exceptions; these reappear in (\ref{nonposex}),
where we consider $p$-groups of non-positive minimum genus.

\section{\bf The main result}\label{proofsec}

In view of the examples worked out in
Sections \ref{smallranksec} and \ref{smallexpsec}, if $G$ runs through
all abelian $p$-groups, there seems to be
a tendency that there are phenomena of `exceptional' and `generic'
cases, where in the `generic' region we have $\mu_0(G)=\mu_e$; 
for an example illustration how this is to be understood see 
Table \ref{regiontbl} (page \pageref{regiontbl}). 
Our main result, to which we proceed in this section, can be seen
as a verification of this observation for a large part of the 
`generic' region. 

We keep the notation introduced in Section \ref{transsec},
in particular let 
$$ G\cong\Z_p^{r_1}\oplus\Z_{p^2}^{r_2}\oplus\cdots\oplus\Z_{p^e}^{r_e} ,$$
where $e\geq 1$, and $r_i\geq 0$ for $1\leq i\leq e-1$, and $r_e\geq 1$.

\subsection{Proposition}\label{minprop}
Suppose that 
$$ \wp(r_{i+1},\dotsc,r_e)\geq p^{e-i}-1 ,$$
for all $0\leq i\leq e-1$ such that $s_{i+1}$ is odd.
Then we have $\mu_0(G)=\mu_e$.

If $s_i>s_{i+1}$ for all $1\leq i\leq e-1$ 
such that $s_{i+1}$ is odd, then the converse also holds.
 
\subsection*{Proof}
By (\ref{minrem}), we have ${\mathrm {min}} ~\gamma (A_e) = \mu_e$
and ${\mathrm {min}} ~\gamma (A_0) = \mu_0$, while for $1\leq i\leq e-1$
we have ${\mathrm {min}} ~\gamma (A_i) \geq \mu_i$. 
Moreover, for $p=2$ and $e'<i\leq e$ we have 
${\mathrm {min}} ~\gamma (A'_i) = \mathrm {min} ~\gamma (A_i)$.
Thus it is sufficient to show that under the assumptions made we
have $\mu_i\geq\mu_e$, for $0\leq i\leq e-1$:

Now $\mu_i\geq\mu_e$ is equivalent to saying
$$ 2\cdot\lfloor\frac{s_{i+1}}{2}\rfloor\cdot p^{e-i}\geq
(p-1) \cdot\wp(s_{i+1},\dotsc,s_e). $$
The right hand side of this inequality being equal to
$$ s_{i+1}p^{e-i}-s_e+\sum_{j=i+1}^{e-1} (s_{j+1}-s_j)p^{e-j}
=s_{i+1}p^{e-i}-1-\wp(r_{i+1},\dotsc,r_e), $$
we thus have $\mu_i\geq\mu_e$ if and only if
$$ (s_{i+1}-2\cdot\lfloor\frac{s_{i+1}}{2}\rfloor)\cdot p^{e-i}
\leq 1+\wp(r_{i+1},\dotsc,r_e) .$$
The latter inequality clearly holds if $s_{i+1}$ is even, 
while if $s_{i+1}$ is odd then it holds if and only if
$\wp(r_{i+1},\dotsc,r_e)\geq p^{e-i}-1$.
This proves the first assertion.

For the second assertion, let $0\leq i\leq e-1$ such that $s_{i+1}$ is odd.
Then for $i\neq 0$ the assumption $s_i-s_{i+1}\geq 1$ implies 
$\epsilon_i=0$, using the notation of (\ref{minrem}), while we have
$\epsilon_0=0$ anyway. Thus we get
$\mu_i={\mathrm {min}} ~\gamma (A_i)\geq\mu_0(G)=\mu_e$,
which by the above observation implies the second assertion. 
\QED

\bigskip 
We are now in a position to prove our main result:

\subsection{Main Theorem}\label{genusthm} 
Let $G$ be a non-trivial abelian $p$-group of shape
$$ G\cong\Z_p^{r_1}\oplus\Z_{p^2}^{r_2}\oplus\cdots\oplus\Z_{p^e}^{r_e} ,$$
such that 
$$ r_i\geq p-1 \quad \textrm{for} \quad 1\leq i\leq e-1, 
\quad \textrm{and} \quad r_e \geq\max\{p-2,1\} .$$

{\bf a)}
Then the reduced minimum and stable upper genera of $G$ are given as
\[ \mu_0(G)= \sigma_0(G) = 
\frac{1}{2}\cdot\left(-1-p^e+\sum_{i=1}^e (p^e-p^{e-i})\cdot r_i \right) .\]

{\bf b)}
Letting $0\leq j \leq e$ be chosen smallest such that
$(r_{j+1},\dotsc,r_e)=(p-1,\dotsc,p-1)$, where $j=e$ refers
to the case $r_e\neq p-1$, the reduced minimum genus $\mu_0(G)$ 
is afforded precisely by the $p$-data
$$ \left(r_1, \dotsc, r_{i-1}, r_i + 1,0,\dotsc,0;\frac{1}{2}(e-i)(p-1)
   \right) ,$$
where $j\leq i\leq e$ is arbitrary for $p$ odd, but
restricted to the cases where $e-i$ is even for $p=2$. 
In particular, $\mu_0(G)$ is always afforded by 
$$ (r_1, \dotsc, r_{e-1}, r_e + 1;0) .$$

\subsection*{Proof} 
{\bf a)} 
By (\ref{minrem}) and (\ref{backrem}) we have 
$$ \frac{1}{2}\cdot\left(-1-p^e+\sum_{i=1}^e (p^e-p^{e-i})\cdot r_i \right) 
=-p^e + {\frac {p-1}{2}} \cdot \wp(s_1,\dotsc,s_e)=\mu_e .$$
Note that $\mu_e\in\frac{1}{2}\Z$, where 
$\mu_e\in\frac{1}{2}\Z\setminus\Z$ if and only if $p=2$ and $s_e$ is odd.
Since $\mu_0(G)\leq\sigma_0(G)$ anyway, it suffices to prove 
$\sigma_0(G)\leq\mu_e$ and $\mu_e\leq\mu_0(G)$:

{\bf i)} 
We first show $\sigma_0(G)\leq\mu_e$: 
By assumption, we have $s_i-s_{i+1}=r_i\geq p-1$ for $1\leq i\leq e-1$, 
that is $||(s_1, \dotsc, s_e)|| \geq p-1$. Hence for any $m\in\N_0$,
by (\ref{mainlinethm}), there is a sequence 
$(a_1, \dotsc, a_e)\in{\mathcal N}(e)$ 
such that $(a_1, \dotsc, a_e)\geq (s_1,\dotsc,s_e)$ and
$\wp(a_1, \dotsc, a_e)=\wp(s_1, \dotsc, s_e)+m$.

Let first $p$ be odd, and $\sigma\in\Z$ such that $\sigma\geq\mu_e$.
Then there are $m\in{\mathbb N}_0$ and
$r\in{\mathbb N}_0$ such that $r<{\frac {p-1}{2}}$ and
$$ \sigma=\mu_e+m\cdot{\frac {p-1}{2}}+r=
-p^e + {\frac {p-1}{2}} \cdot (\wp(s_1, \dotsc, s_e)+m)+r .$$
Let $(a_1, \dotsc, a_e)$ as above such that 
$\wp(a_1, \dotsc, a_e)=\wp(s_1, \dotsc, s_e)+m$, and $a_{e+1}:=2r$,
then $a_e-a_{e+1}\geq (r_e+1)-2\cdot{\frac {p-3}{2}}\geq 2$ 
implies $(a_1,\dotsc,a_{e+1})\in A_e$. Since 
$\gamma(a_1,\dotsc,a_{e+1})=
-p^e + r + {\frac {p-1}{2}} \cdot \wp(a_1, \dotsc, a_e) = \sigma$,
from (\ref{talurem}) we get $\sigma\in\textrm{sp}_0(G)$.

Let now $p=2$, and $\sigma\in\frac{1}{2}\Z$ such that $\sigma\geq\mu_e$.
Let $m:=2(\sigma-\mu_e)\in\N_0$. Let $(a_1,\dotsc,a_e)$ be as above such that 
$\wp(a_1, \dotsc, a_e)=\wp(s_1, \dotsc, s_e)+m$, and $a_{e+1}:=0$, then 
$a_e-a_{e+1}\geq r_e+1\geq 2$ implies $(a_1,\dotsc,a_{e+1})\in A_e$. 
Since $\gamma(a_1,\dotsc,a_{e+1})=
-2^e + {\frac {1}{2}} \cdot \wp(a_1, \dotsc, a_e) = \sigma$.
Thus, if $e'=e$ 
from (\ref{talurem}) we get $\sigma\in\textrm{sp}_0(G)$.

If $e'<e$, then we have $e'=e-1$ and $s_e=2$, and hence
$\gamma(A(G))=\gamma(A_0)\cup \gamma(A_1) \cup \cdots \cup 
\gamma(A_{e-1}) \cup \gamma(A'_e)\subseteq\Z$.
Since $\mu_e={\mathrm {min}} ~\gamma (A'_e)$ we may assume that 
$\sigma\in\Z$, thus $m:=2(\sigma-\mu_e)\in\N_0$ is even. Hence we get
$$ a_e\equiv\wp(a_1, \dotsc, a_e)=\wp(s_1, \dotsc, s_e)+m 
\equiv s_e+m \equiv 0\pmod{2} ,$$
implying that $(a_1,\dotsc,a_{e+1})\in A'_e$, and 
from (\ref{talurem}) we get $\sigma\in\textrm{sp}_0(G)$.

{\bf ii)} 
We show $\mu_e \leq \mu_0(G)$:
Since $s_i-s_{i+1}=r_i\geq 1$ for all $1\leq i\leq e-1$, 
by (\ref{minprop}) we have to show 
$\wp(r_{i+1},\dotsc,r_e) \geq  p^{e-i} -1$,
for all $0\leq i\leq e-1$ such that $s_{i+1}$ is odd.

For $p$ odd we have $r_j\geq p-1$ for $1\leq j\leq e-1$, and $r_e\geq p-2$,
where $\sum_{j=i+1}^e r_j=s_{i+1}-1$ being even implies that
$(r_{i+1},\dotsc,r_{e-1,}r_e)\neq (p-1,\dotsc,p-1,p-2)$. Thus
$$ \wp(r_{i+1},\dotsc,r_e)>-1+(p-1)\cdot\sum_{j=i+1}^e p^{e-j}=p^{e-i}-2 .$$
For $p=2$ we have $r_j\geq 1$ for $1\leq j\leq e$, directly yielding
$$ \wp(r_{i+1},\dotsc,r_e)=\sum_{j=i+1}^e r_j\cdot 2^{e-j} \geq 
\sum_{j=i+1}^e 2^{e-j} = 2^{e-i}-1 .$$

{\bf b)} 
We determine when $\mu_0(G)$ is attained:
By (\ref{minrem}), ${\mathrm {min}} ~\gamma (A_e) = \mu_e$
is attained precisely for $(s_1,\dotsc,s_e,0)$, corresponding
to the $p$-datum $(r_1, \dotsc, r_{e-1}, r_e + 1;0)$.

Now, for $0\leq i\leq e-1$, 
by the proof of (\ref{minprop}) we have $\mu_i\geq\mu_e$.
Moreover, replacing inequalities by equalities in the proof of (\ref{minprop}) 
shows that $\mu_i=\mu_e$ is equivalent to $s_{i+1}$ being odd and
$\wp(r_{i+1},\dotsc,r_e)=p^{e-i}-1$.
Since $(r_{i+1},\dotsc,r_{e-1},r_e)\geq (p-1,\dotsc,p-1,\max\{p-2,1\})$,
the latter equality holds if and only if
$(r_{i+1},\dotsc,r_e)=(p-1,\dotsc,p-1)$. Since in this case 
$s_{i+1}-1=\sum_{j=i+1}^e r_j=(e-i)(p-1)$, we have $s_{i+1}$ odd 
if and only if $p$ is odd or $e-i$ is even.
Hence we conclude, by (\ref{minrem}) again, that in these cases
${\mathrm {min}} ~\gamma (A_i) = \mu_i$ is attained precisely for 
$$ (s_1, \dotsc, s_i, 2\cdot \lfloor {\frac {s_{i+1}}{2}} \rfloor, \dotsc, 
2\cdot \lfloor {\frac {s_{i+1}}{2}} \rfloor)
= (s_1, \dotsc, s_i, s_{i+1}-1, \dotsc, s_{i+1}-1) ,$$
corresponding to the $p$-datum, using the notation of (\ref{backrem}), 
$$ \un{x}^i
= (r_1, \dotsc, r_{i-1}, r_i + 1,0,\dotsc,0;\frac{1}{2}(e-i)(p-1)) .$$
Note that we have $\mathcal I(G)=\{0,\dotsc,e\}$ for $p$ odd,
while for $p=2$ we at least get 
$\{0\}\cup\{e-2\cdot\lfloor\frac{e-j}{2}\rfloor,\dotsc,e-2,e\}
\subseteq\mathcal I(G)$,
hence the indices $0\leq i\leq e$ affording $\mu_0(G)$ are indeed elements of
the index set $\mathcal I(G)$, in accordance with (\ref{minthm}).
\QED

\subsection{Example}
{\bf i)}
For $p$ odd and $(r_1,\dotsc,r_{e-1},r_e)=(p-1,\dotsc,p-1,p-2)$,
that is the extremal case, we get, recovering \cite[Cor.3.7]{tal},
$$ \mu_0(G)=\sigma_0(G)=
\frac{1}{2}\cdot(((e(p-1)-3)\cdot p^e+1) .$$

{\bf ii)}
For $p$ arbitrary and $(r_1,\dotsc,r_{e-1},r_e)=(p-1,\dotsc,p-1,p-1)$ we get
$$ \mu_0(G)=\sigma_0(G)=\frac{1}{2}\cdot(e(p-1)-1)\cdot p^e ,$$
which for $p=2$ specializes to $\mu_0(G)=\frac{e-2}{2}\cdot 2^e$. 

\bigskip

As an immediate consequence of (\ref{genusthm}), invoking
Kulkarni's Theorem (\ref{kulkarnithm}), we are able to describe
the complete (reduced) spectrum of the groups in question:

\subsection{Corollary}\label{genuscor} 
{\bf a)}
The reduced spectrum of $G$ is given as
$$ {\mathrm{sp}}_0(G)=
\left\{\begin{array}{ll}
\mu_0(G)+{\mathbb N}_0, & \textrm{if }p\textrm{ odd  or }r_e=1, \\
\mu_0(G)+\frac{1}{2}{\mathbb N}_0, & \textrm{if }p=2\textrm{ and }r_e\geq 2. \\
\end{array}\right. $$

{\bf b)}
Letting $\delta=\delta(G):=\sum_{i=1}^e (ir_i-1)$ 
be the cyclic deficiency of $G$, then 
the minimum genus and the spectrum of $G$ are given as 
$\mu(G)=1+p^\delta\cdot\mu_0(G)$ and
$$ {\mathrm{sp}}(G)=
\left\{\begin{array}{ll}
1+p^\delta\cdot\mu_0(G)+p^\delta\cdot{\mathbb N}_0, 
& \textrm{if }p\textrm{ odd or }r_e=1, \\
1+2^\delta\cdot\mu_0(G)+2^{\delta-1}\cdot{\mathbb N}_0, 
& \textrm{if }p=2\textrm{ and }r_e\geq 2. \\
\end{array}\right. $$

\bigskip
Moreover, for certain suitable co-finite sets of positive integers we are 
conversely able to provide abelian $p$-groups having the specified set 
as their reduced spectrum:

\subsection{Theorem}\label{genspecthm}
Let $p$ be a prime, let $e \geq 1$, and let $m\in\N$ such that  
$$ m\geq \left\{
\begin{array}{ll} 
(2e-1)p^e-2\cdot\frac{p^e-1}{p-1}+1,
& \textrm{if }p\textrm{ odd}, \\
(e-1)\cdot 2^{e+1}+2, & \textrm{if }p=2. \rule{0em}{1.2em} \\
\end{array}\right. $$
Then there is a group $G$ of exponent $p^e$ such that
$\mu_0(G) = -p^e + {\frac {p-1}{2}} \cdot m$ and 
$$ {\mathrm {sp}}_0(G)=\left\{
\begin{array}{ll} 
\mu_0(G)+{\mathbb N}_0, & \textrm{if }p\textrm{ odd or }m\textrm{ even}, \\
\mu_0(G)+\frac{1}{2}{\mathbb N}_0, & 
         \textrm{if }p=2\textrm{ and }m\textrm{ odd} .\\
\end{array}\right.  $$

\subsection*{Proof}
We consider the sequence $(a_1,\dotsc,a_e)\in{\mathcal N}$ 
given by $a_e:=\max\{p-1,2\}$,
and $a_{e-i}:=a_e+i\cdot 2(p-1)$ for $1\leq i\leq e-1$. 

{\bf i)}
We first show that the lower bound for $m$ given above 
coincides with $\wp(a_1,\dotsc,a_e)$:

To this end, we first observe that
$s_e(p):=\sum_{i=1}^e ip^i
=\frac{p}{p-1}\cdot(ep^e-\sum_{i=0}^{e-1}p^i)$,
which in turn is seen by induction: 
This formula being correct for $e=1$, we get
$s_{e+1}(p)=(e+1)p^{e+1}+s_e(p)
=\frac{p}{p-1}\cdot\left((e+1)(p-1)p^e+ep^e-\sum_{i=0}^{e-1}p^i\right)
=\frac{p}{p-1}\cdot\left((e+1)p^{e+1}-\sum_{i=0}^e p^i\right)$.
In particular, for $p=2$ we get 
$s_e(2)=(e-1)\cdot 2^{e+1}+2$.

Now, for $p$ odd we have
$$ \wp(a_1,\dotsc,a_e)=(p-1)\cdot\sum_{i=1}^e(2(e-i)+1)p^{e-i}
=\frac{2(p-1)}{p}\cdot\sum_{i=1}^e ip^i
-(p-1)\cdot\sum_{i=0}^{e-1} p^i ,$$
which using the above expression for $s_e(p)$ can be rewritten as
$$\wp(a_1,\dotsc,a_e)
=2\cdot\left(ep^e-\sum_{i=0}^{e-1}p^i\right)-p^e+1 
=(2e-1)p^e+1-2\cdot\sum_{i=0}^{e-1}p^i .$$ 
For $p=2$ we get
$\wp(a_1,\dotsc,a_e)=2\cdot\sum_{i=1}^e(e-i+1)\cdot 2^{e-i}
=\sum_{i=1}^e i\cdot 2^i=s_e(2)$.

{\bf ii)}
The strategy of proof now is reminiscent of the proof of (\ref{mainlinethm}): 
Given $m \geq \wp(a_1,\dotsc,a_e)$, then we write
$m-\wp(a_1,\dotsc,a_e)$ in a partial $p$-adic expansion as
$m-\wp(a_1,\dotsc,a_e) = \sum_{i=1}^e b_ip^{e-i}$,
where $b_i\geq 0$ such that $b_2,\dotsc,b_e<p$, 
but $b_1$ might be arbitrarily large. 
Hence letting $s_i:=a_i+b_i$ for $1\leq i\leq e$,
we have $m=\sum_{i=1}^e s_ip^{e-i}$.
Thus for $1\leq i\leq e-1$ we get 
$$r_i:=s_i-s_{i+1} = 2(p-1) + (b_i - b_{i+1}) \geq p-1, $$
and $r_e:=s_e-1\geq a_e-1=\max\{p-2,1\}$. 
Hence, by (\ref{genusthm}), for the abelian group of shape 
$G \cong {\mathbb Z}_p^{r_1} \oplus {\mathbb Z}_{p^2}^{r_2} 
\oplus \dotsc \oplus {\mathbb Z}_{p^e}^{r_e}$
we have 
$$ \sigma_0(G)=\mu_0(G)=\mu_e
=-p^e + {\frac {p-1}{2}} \cdot \wp(s_1,\dotsc,s_e)
=-p^e + {\frac {p-1}{2}} \cdot m .$$
Moreover, for $p=2$ we have $a_e=2$, and thus 
if $m$ is even we get $b_e=0$ and hence $r_e=1$,
while if $m$ is odd we get $b_e=1$ and hence $r_e=2$. 
Thus the statement on ${\mathrm {sp}}_0(G)$ follows from (\ref{genuscor}).
\QED

\section{\bf Talu's Conjecture}\label{applsec}

In general, we might wonder which invariants of a non-trivial
abelian $p$-group $G$ are determined by its spectrum. Given the latter,
this determines the Kulkarni invariant $N=N(G)$,
and hence the cyclic deficiency $\delta=\delta(G)=\log_p(N)$ 
is known as well whenever $p$ is odd, while 
$\delta\in\{\log_p(N),1+\log_p(N)\}$ for $p=2$. Thus the spectrum also
determines the reduced minimum and stable upper genera whenever $p$ is odd,
while the latter are known up to a factor of $2$ for $p=2$.

In this spirit, Talu's Conjecture says that, if $p$ is odd, then even
the isomorphism type of $G$ is determined by its spectrum. 
We are tempted to include the case $p=2$ as well by expecting this to
hold true up to finitely many finite sets of exceptions; we cannot
possibly expect more, for example in view of the sets of groups
$\{\Z_2,\Z_4,\Z_2^2,\Z_8\}$ and $\{\Z_2\oplus\Z_4,\Z_2^3,\Z_2\oplus\Z_8\}$
discussed in (\ref{z4p2spec}).

As for evidence, restricting to certain classes of abelian
$p$-group, Talu's Conjecture (including the case $p=2$) holds 
within the class of cyclic $p$-groups with the only exception
of $\{\Z_2,\Z_4,\Z_8\}$, see (\ref{cycex});
within the class of elementary abelian $p$-groups with the only exception
of $\{\Z_2,\Z_2^2\}$, see (\ref{e1ex});
and within the class of $p$-groups of exponent $p^2$, see (\ref{e2uniq}).  
We proceed to prove a further positive result:

\subsection{A finiteness result}\label{finrem}
We show that, as long as we stick to groups fulfilling the assumptions
of (\ref{genusthm}), given the spectrum of $G$ there are only finitely
many groups having the same spectrum, up to isomorphism. Actually, 
just keeping the reduced minimum genus fixed leaves only finitely
possibilities:

Note first that the only admissible cyclic groups are $\Z_2$ and $\Z_3$, 
hence we may assume that the groups 
we are looking for are non-cyclic, that is have an associated
sequence $(s_1,\dotsc,s_e)\neq (2,\dotsc,2)$.
We now show that, given any $m\geq 0$, there are only finitely many
$e\geq 1$ and sequences $s_1 \geq \cdots \geq s_e \geq 2$, where
$s_1\geq 3$, such that
$$ \mu_e=-p^e + \frac{p-1}{2}\cdot\wp(s_1, \dotsc, s_e)\leq m .$$

This is seen as follows: The above inequality is equivalent to
$$\wp(s_1-2, \dotsc, s_e-2) =
 \wp(s_1, \dotsc, s_e) -2\cdot\frac{p^e-1}{p-1} \leq \frac {2(m+1)}{p-1} .$$
This implies $(s_1-2)\cdot p^{e-1}\leq \frac {2(m+1)}{p-1}$, 
hence since $s_1\geq 3$ we infer that $e$ is bounded. Fixing $e$, we get
$(s_i-2)\cdot p^{e-i}\leq \frac {2(m+1)}{p-1}$, bounding $s_i$ as well,
for $1\leq i\leq e$.
\QED

\bigskip

In view of this, there necessarily are groups fulfilling 
the assumptions of (\ref{genusthm}) whose reduced minimum genus exceeds 
any given bound. Hence the point of (\ref{genspecthm}) is to add some
precision to this observation.
But here positive results come to an end: 

In (\ref{e3ex}) and (\ref{counterex}),
we are going to construct counterexamples to
Talu's Conjecture (both for $p$ odd and $p=2$), consisting
of pairs of groups having the same order exponent, and
pairs where these invariants are different, respectively.
Even worse, by the results in (\ref{e3ex}), there cannot be
an absolute bound on the cardinality of a set of abelian $p$-groups
having the same spectrum, not even if we restrict 
to groups having the same order and exponent.

\subsection{Counterexamples with fixed exponent}\label{e3ex}
We construct non-isomorphic abelian groups $G$ and $\tilde G$
having the same order, exponent, and spectrum,
thus in particular having the same Kulkarni invariant, cyclic deficiency,
minimum genus and reduced minimum genus.

In view of the results in (\ref{e1ex}) and (\ref{e2uniq}), we let $e:=3$,
and look at groups 
$G \cong {\mathbb Z}_p^{r_1} \oplus {\mathbb Z}_{p^2}^{r_2} 
\oplus {\mathbb Z}_{p^3}^{r_3}$
and 
$\tilde G \cong {\mathbb Z}_p^{\tilde r_1} \oplus 
{\mathbb Z}_{p^2}^{\tilde r_2} \oplus {\mathbb Z}_{p^3}^{\tilde r_3}$
of exponent $p^3$ fulfilling the assumptions of (\ref{genusthm}),
that is coming from sequences $\un{r}=(r_1,r_2,r_3)$ and
$\un{\tilde r}=(\tilde r_1,\tilde r_2,\tilde r_3)$ such that 
$r_1,r_2,\tilde r_1,\tilde r_2 \geq p-1$ and $r_3,\tilde r_3\geq\max\{p-2,1\}$. 
Then, by (\ref{genuscor}), the groups $G$ and $\tilde G$ are as desired
if and only if they are non-isomorphic such that $|G|=|\tilde G|$ and 
$\mu_0(G)=\mu_0(\tilde G)$, and in case $p=2$ we have 
$r_3=1$ if and only if $\tilde r_3=1$.

Now $|G|=|\tilde G|$ translates into
$$ r_1+2r_2+3r_3=\log_p(|G|)=\log_p(|\tilde G|)=
\tilde r_1+2\tilde r_2+3\tilde r_3 ,$$ 
and $\mu_0(G)=\mu_0(\tilde G)$ translates into
$$ \sum_{i=1}^3 (p^3-p^{3-i})\cdot r_i = 
\sum_{i=1}^3 (p^3-p^{3-i})\cdot \tilde r_i. $$  
Hence we conclude that we have $|G|=|\tilde G|$ and 
$\mu_0(G)=\mu_0(\tilde G)$ if and only if
$\un{\tilde r}-\un{r}\in\Z^3$ is an element of the row kernel of the matrix 
$$ P:=\begin{bmatrix} 1&p^3-p^2 \\ 2&p^3-p\\ 3&p^3-1\\
\end{bmatrix} \in \Z^{3\times 2}\subseteq {\mathbb Q}^{3\times 2} .$$
Now $P$ has $\mathbb Q$-rank $2$, and its row kernel is given as 
$\ker(P)=\langle\un{\rho}\rangle_{\mathbb Q}$, where 
$$ \un{\rho}:=(p+2,-2p-1,p)\in\Z^3 .$$
Since $\gcd(p+2,-2p-1,p)=1$ we conclude that
$\ker(P)\cap \Z^3=\langle\un{\rho}\rangle_\Z$. 

In conclusion, we have $|G|=|\tilde G|$ and 
$\mu_0(G)=\mu_0(\tilde G)$ if and only if
$\un{\tilde r}=\un{r}+k\cdot\un{\rho}$ for some $k\in\Z$,
where $G$ and $\tilde G$ are non-isomorphic if and only if $k\neq 0$.
Thus this provides a complete picture of the counterexamples to 
Talu's Conjecture in the realm of abelian groups
of exponent $p^3$ fulfilling the assumptions of (\ref{genusthm}).
In particular, for any $l\in\N$ there is a set of isomorphism types 
of cardinality at least $l+1$ consisting of groups having the same 
order and reduced minimum genus: 
Given $r_1\geq p-1$ and $r_3\geq p-2$, such that $r_3\geq 2$ for $p=2$,
and letting $r_2:=(p-1)+l\cdot(2p+1)$, all the sequences
$\un{r}+k\cdot\un{\rho}$, where $0\leq k\leq l$, give rise to 
groups as desired.
The smallest counterexamples, in terms of group order, are given by
choosing $\un{r}$ as small as possible for the case $l=1$:

{\bf i)}
For $p$ odd this yields
$$ \un{r}=(p-1,3p,p-2) \quad \textrm{and} \quad 
\un{\tilde r}:=\un{r}+\un{\rho} =(2p+1,p-1,2p-2) ,$$
giving rise to groups such that
$$ |G|=|\tilde G|=p^{10p-7} \quad \textrm{and} \quad
\mu_0(G)=\mu_0(\tilde G)=\frac{1}{2}\cdot (5p^4-5p^3-2p^2-p+1) .$$
Hence in particular for $p=3$ we get $\un{r}=(2,9,1)$ and
$\un{\tilde r}=(7,2,4)$, giving rise to groups such that
$|G|=|\tilde G|=3^{23}$ and $\mu_0(G)=\mu_0(\tilde G)=125$.

{\bf ii)}
In order to cover the case $p=2$ as well, for $p$ arbitrary we may let 
$$ \un{r}=(p-1,3p,p) \quad \textrm{and} \quad \un{\tilde r}=(2p+1,p-1,2p) ,$$
giving rise to groups such that
$$ |G|=|\tilde G|=p^{10p-1} \quad \textrm{and} \quad
\mu_0(G)=\mu_0(\tilde G)=\frac{1}{2}\cdot (5p^4-3p^3-2p^2-p-1) .$$ 
Hence in particular for $p=2$ we get $\un{r}=(1,6,2)$ and
$\un{\tilde r}=(5,1,4)$, giving rise to groups such that
$|G|=|\tilde G|=2^{19}$ and $\mu_0(G)=\mu_0(\tilde G)=\frac{45}{2}$.

\subsection{Counterexamples with varying exponent}\label{counterex}
We construct non-isomorphic abelian groups $G$ and $\tilde G$ just
having the same spectrum, thus in particular having the same Kulkarni
invariant and minimum genus; hence for $p$ odd also having the same 
cyclic deficiency and reduced minimum genus. We might wonder whether 
in this situation, possibly further assuming that $G$ and $\tilde G$
belong to the class of groups described in (\ref{genusthm}), the
groups necessarily have the same exponent, or equivalently the same 
order whenever $p$ odd; if this was the case then the examples in 
(\ref{e3ex}) would be the typical or even the only ones.

We look at groups afforded by sequences $\un{r}=(r_1,\dotsc,r_e)$
and $\un{\tilde r}=(\tilde r_1,\dotsc,\tilde r_{\tilde e})$, where
$e,\tilde e\geq 1$, fulfilling the assumptions of (\ref{genusthm}), 
that is $r_i\geq p-1$ for $1\leq i\leq e-1$, and
$\tilde r_i\geq p-1$ for $1\leq i\leq \tilde e-1$, as well as 
$r_e,\tilde r_{\tilde e}\geq\max\{p-2,1\}$. 
We are going to present a series of counterexamples to Talu's Conjecture 
fulfilling $e\neq\tilde e$, where this subsection deals with 
the case $p$ odd, while the case $p=2$ is treated in (\ref{counterex2}).
But before doing so, we would like to indicate the
heuristics we have used to find them:

Let $\delta\geq -2e+\frac{e(e+1)}{2}\cdot (p-1)$ whenever $p$ is odd,
and $\delta\geq \frac{e(e-1)}{2}$ for $p=2$,
in each case the lower bound being the cyclic deficiency associated with
the smallest admissible sequence $(p-1,\dotsc,p-1,\max\{p-2,1\})$;
note that smaller values of $\delta$ are not achieved at all.
We now aim at varying $\un{r}$ within the set of admissible sequences,
such that $\log_p(|G|)=\delta+e=\sum_{i=1}^e ir_i$ is kept fixed, but
$$ 2\mu_e+1=-p^e+\sum_{i=1}^e (p^e-p^{e-i})\cdot r_i 
=-p^e+\sum_{i=1}^e \frac{p^e-p^{e-i}}{i}\cdot ir_i $$
is maximized and minimized, respectively.

To this end, we observe that the arithmetic mean of the first $i$
entries of the sequence $(p^{e-1},\dotsc,p,1)$ is given as
$\frac{1}{i}\cdot\sum_{j=e-i}^{e-1}p^j
=\frac{1}{i}\cdot\frac{p^e-p^{e-i}}{p-1}$,
for $1\leq i\leq e$, hence the sequence
$(\frac{p^e-p^{e-1}}{1},\frac{p^e-p^{e-2}}{2},\dotsc,\frac{p^e-1}{e})$
is strictly decreasing. Thus $2\mu_e+1$ becomes largest (respectively
smallest) by choosing the last (respectively first) $e-1$ entries of 
$\un{r}$ as small as possible, and adjusting the first (respectively last)
entry such that $\un{r}$ has cyclic deficiency $\delta$ associated with it.

For the remainder of this subsection let $p$ be odd. Then maximizing yields
$2\mu_e+1\leq 2\mu_e(a,p-1,\dotsc,p-1,p-2)+1$, where 
$$ a:=\delta+2e-\frac{(e+2)(e-1)}{2}\cdot (p-1) .$$ 
Note that by the choice of $\delta$ we conclude that $a\geq p-1$, hence
the right hand side of the above inequality is achieved. 
By a straightforward computation we get 
$$ \begin{array}{rcl}
2\mu_e+1 & \leq & \quad
\left(\delta+\frac{(e-1)(e+6)}{2}-\frac{e(e-1)(p-1)}{2}\right)\cdot p^e \\
&& -\left(\delta+\frac{e(e+5)}{2}\right)\cdot p^{e-1}+2\rule{0em}{1.7em} \\
\end{array} $$

Similarly, minimizing yields
$2\mu_e+1\geq 2\mu_e(p-1,p-1,\dotsc,p-1,b)+1$, where 
$$ b:=\frac{\delta}{e}-\frac{e-1}{2}\cdot (p-1)+1 .$$ 
Note that here $b$ in general is not integral, so that the right hand side
of the above inequality might not be achieved; 
it is possible to determine explicitly the sequence giving rise
to the actual minimum of $2\mu_e+1$, but this will not be needed.
By a straightforward computation we get 
$$ 2\mu_e+1 \geq
\left(\frac{\delta}{e}+\frac{(e-1)(p-1)}{2}-1\right)\cdot p^e
+\frac{(e+1)(p-1)}{2}-\frac{\delta}{e} .$$

Hence we have to ensure that the above upper bound for 
$2\mu_{\tilde e}+1$, applied to some $1\leq \tilde e<e$, 
is at least as large as the lower bound for $2\mu_e+1$. 
Viewing the upper and lower bounds as linear functions in $\delta$,
in order to have an unbounded range of candidates $\delta$ to check,
the slope of the upper bound function should exceed the slope
of the lower bound function. This yields
$$ (p-1)p^{\tilde e-1} \geq \frac{p^e-1}{e} ,$$
in other words
$$ e \geq 
\sum_{i=1}^e p^{(e-i)-(\tilde e-1)} =
\sum_{i=0}^{e-\tilde e} p^i +
\sum_{i=1}^{\tilde e-1} p^{-i} =
\frac{p^{e-\tilde e+1}-1}{p-1}+\sum_{i=1}^{\tilde e-1} p^{-i} ,$$
implying 
$$ e\geq \frac{p^{e-\tilde e+1}-1}{p-1}+1 .$$

Thus we are led to consider the case $\tilde e=e-1$, where 
the smallest possible choices are $e:=p+2$ and $\tilde e:=p+1$.
This yields the following specific examples: Let
$$ \un{r}:=(p-1,\dotsc,p-1,p,p^3+p^2-2) ,$$
thus having $p$ consecutive entries $p-1$, and for $p\geq 5$ let
$$ \tilde{\un{r}}:=(p^4+3p^3+2p^2-p-1,p-1,\dotsc,p-1,p,p,p-1,p-2) ,$$ 
thus having $p-4$ consecutive entries $p-1$, while for $p=3$ let
$$ \tilde{\un{r}}:=(p^4+3p^3+2p^2-p,p,p-1,p-2)|_{p=3}=(177,3,2,1) ;$$
a few explicit cases are given in Table \ref{counterextbl}.
Then, by a straightforward computation, we indeed have 
$$ \delta=\tilde\delta=p^4+\frac{7}{2}p^3+3p^2-\frac{5}{2}p-6 ,$$
and
$$ \mu_e(\un{r}) = \mu_{\tilde e}(\tilde{\un{r}}) = \frac{1}{2}\cdot
   \left( (p^3+2p^2-4)\cdot p^{p+2}-p^3-p^2+1\right) .$$ 
We remark that, had we carried out the analysis on minimizing 
$2\mu_e+1$, we would have found $\un{r}$ as the minimizing sequence
associated with $\delta$.
Thus $\un{r}$ and $\tilde{\un{r}}$ 
give rise to groups $G$ and $\tilde G$, respectively,
by (\ref{genuscor}) having the same spectrum, 
but having distinct exponents $p^{p+2}$ and $p^{p+1}$, respectively.
 
Actually, the above series has been found by running an explicit search 
for odd $p\leq 11$, using the computer algebra system {\sf GAP} \cite{GAP},
and observing the pattern arising. We suspect that these in general are the 
counterexamples to Talu's Conjecture with smallest possible cyclic 
deficiency $\delta$ for groups of exponents $p^{p+2}$ and $p^{p+1}$, 
respectively; but we have not attempted to prove this in general,
and only checked it explicitly for $3\leq p\leq 23$ using {\sf GAP}.

The above analysis also implies that counterexamples consisting of 
groups of exponent $p^e$ and $p^{\tilde e}$, respectively, such that 
$\tilde e <e\leq p+1$ can possibly exist only for finitely many 
values of $\delta$. Actually, we suspect that counterexamples such that 
$\tilde e\leq p$ do not exist at all; but we have not thoroughly
investigated into this, and only made a few unsuccessful explicit
searches for $3\leq p\leq 23$ using {\sf GAP}.

\begin{table}[t]\caption{Counterexamples with varying exponent.}
\label{counterextbl}
\vspace*{-1.5em}

$$ \begin{array}{|r|cc|}
\hline
p & \un{r} & \tilde{\un{r}} \\
\hline
\hline
3 & (2, 2, 2, 3, 34) & (177, 3, 2, 1) \\
5 & (4,\dotsc,4,5,148) & (1044, 4, 5, 5, 4, 3) \\
7 & (6,\dotsc,6,7,390) & (3520, 6, 6, 6, 7, 7, 6, 5 ) \\ 
11 & (10, \dotsc, 10, 11, 1450 ) & (18864,10,\dotsc,10,11,11,10,9 ) \\
13 & (12, \dotsc, 12, 13, 2364 ) & (35476,12,\dotsc,12,13,13,12,11 ) \\
17 & (16, \dots, 16, 17, 5200) & (98820,16,\dotsc, 16, 17, 17, 16, 15) \\
\hline
\end{array} $$

$$ \begin{array}{|r|r|r|r|}
\hline
p & e & \delta & \mu_e \\
\hline
\hline
3 & 5 & 189 & 4964 \\
5 & 7 & 1119 & 6679613 \\
7 & 9 & 3725 & 8817262934 \\
11 & 13 & 19629 & 27083067676913144 \\
13 & 15 & 36719 & 64775747609331851801 \\
17 & 19 & 101535 & 655895227302212659718161655 \\
\hline
\end{array} $$

\vspace*{0.5em}
\hrulefill
\end{table}

\subsection{Counterexamples with varying exponent for $p=2$}\label{counterex2}
We keep the setting of (\ref{counterex}), but let now $p=2$. 
Since our approach involves sequences $\un{r}$ such that $r_e\geq 2$,
for $\tilde{\un{r}}$ we distinguish the cases $\tilde r_{\tilde e}\geq 2$
and $\tilde r_{\tilde e}=1$:

{\bf i)}
Let first $\tilde r_{\tilde e}\geq 2$. Then, by (\ref{genuscor}),
the groups $G$ and $\tilde G$ associated with these sequences 
have the same spectrum if and only if they
have the same cyclic deficiency and reduced minimum genus.
Thus a similar analysis as the one in the odd prime case yields
$2\mu_e+1\leq 2\mu_e(a,1,\dotsc,1)+1$, where 
$a:=\delta-\frac{(e-2)(e+1)}{2}$, hence we get
$$ 2\mu_e+1 \leq 
\left(\delta-\frac{(e-2)(e-3)}{2}-1\right)\cdot 2^{e-1}+1 .$$
Similarly, we get $2\mu_e+1\geq 2\mu_e(1,\dotsc,1,b)+1$, where 
$b:=\frac{\delta}{e}-\frac{e-3}{2}$, yielding
$$ 2\mu_e+1 \geq
\left(\frac{\delta}{e}+\frac{e-3}{2}\right)\cdot 2^e
+\frac{e+1}{2}-\frac{\delta}{e} .$$
Again comparing slopes with respect to $\delta$ of the upper and lower
bound functions yields $2^{\tilde e-1}\geq \frac{2^e-1}{e}$, which is 
the same formula as in the odd prime case, specialized to $p=2$.
Hence here we obtain the condition $e\geq 2^{e-\tilde e+1}$.
Moreover, it turns out that for $1\leq \tilde e<e\leq 3$ and
any $\delta\geq 0$ the upper bound for $2\mu_{\tilde e}+1$ is
smaller than the lower bound for $2\mu_e+1$, excluding these 
choices of $\tilde e<e$. Hence we are led to consider the case 
$\tilde e=e-1$, with smallest possible choices $e:=4$ and $\tilde e=3$:

An explicit search using {\sf GAP} yields the smallest counterexamples, 
with respect to cyclic deficiency $\delta$, as
$$ \un{r}:=(1, 1, 1, 18) 
\quad\textrm{and}\quad
\tilde{\un{r}}:=(69, 1, 2).$$
Then we get $\delta=\tilde\delta=74$ and 
$\mu_e(\un{r}) = \mu_{\tilde e}(\tilde{\un{r}}) = \frac{287}{2}$,
where again we remark that $\un{r}$ is the minimizing sequence
associated with $\delta$.
Thus $\un{r}$ and $\tilde{\un{r}}$ 
give rise to groups $G$ and $\tilde G$, respectively,
by (\ref{genuscor}) having the same spectrum, and both 
fulfilling the `$e'=e$' property, 
but having distinct exponents $16$ and $8$, respectively.

{\bf ii)}
Let now $\tilde r_{\tilde e}=1$. Then, by (\ref{genuscor}) the groups 
$G$ and $\tilde G$ associated with the sequences $\un{r}$
and $\tilde{\un{r}}$ have the same spectrum if and only if for the
associated cyclic deficiency and reduced minimum genus we have 
$$ \tilde\delta=\delta-1
\quad\textrm{and}\quad
\mu_{\tilde e}(\tilde{\un{r}})=2\mu_e(\un{r}) .$$
Considering again the slopes with respect to $\delta$ 
of the upper and lower bound functions, from 
$2\mu_{\tilde e}(\tilde{\un{r}})+1= 4\mu_e(\un{r})+1 
=2\cdot (2\mu_e(\un{r})+1)-1$
we this time get $2^{\tilde e-1}\geq 2\cdot \frac{2^e-1}{e}$,
implying $e\geq 2^{e-\tilde e+2}$,
thus leading us to consider the case $\tilde e=e-1$
with smallest possible choices $e:=8$ and $\tilde e=7$:

An explicit search using {\sf GAP} yields the smallest counterexamples,
with respect to cyclic deficiency $\delta$, as
$$ \un{r}:=(1, 1, 1, 1, 1, 1, 1, 1025) 
\quad\textrm{and}\quad
\tilde{\un{r}}:=(8199, 1, 1, 1, 1, 1, 1).$$
Then we get $\delta=8220=\tilde\delta+1$ and 
$\mu_e(\un{r}) = 131328 =\frac{1}{2}\cdot \mu_{\tilde e}(\tilde{\un{r}})$,
where again we remark that $\un{r}$ is the minimizing sequence 
associated with $\delta$, and $\tilde{\un{r}}$ is the maximizing
sequence associated with $\tilde\delta$. Thus $\un{r}$ and $\tilde{\un{r}}$ 
give rise to groups $G$ and $\tilde G$, respectively,
by (\ref{genuscor}) having the same spectrum, precisely one of them 
fulfilling the `$e'=e$' property,
and having distinct exponents $256$ and $128$, respectively.
Moreover, although we have not thoroughly investigated into this,
unsuccessful explicit searches using {\sf GAP} lead us to suspect
that such counterexamples with $\tilde e \leq 6$ do not exist.
\QED

\bigskip 

Finally, we remark that the above approach can also be used
to find counterexamples fulfilling $\tilde e=e$: 
Actually, by (\ref{e1ex}) and (\ref{e2uniq}), there cannot
be counterexamples for $1\leq \tilde e = e\leq 2$,
except the groups $\{\Z_2,\Z_2^2\}$; note that the latter indeed
is a single counterexample, for $\delta=1$,
while our approach is aiming at finding $\tilde e\leq e$
allowing for an infinite range of candidates $\delta$.
Moreover, it turns out that for $\tilde e=e=3$ and
any $\delta\geq 0$ the upper bound for $2\mu_e+1$ is
smaller than the lower bound for $2\cdot(2\mu_e+1)-1$, 
excluding this case. Hence we are led to consider the case $e:=4$:

An explicit search using {\sf GAP} yields the smallest examples,
with respect to cyclic deficiency $\delta$, as
$$ \un{r}:=(1, 1, 1, 21) 
\quad\textrm{and}\quad
\tilde{\un{r}}:=(80, 1, 1, 1).$$
Then we get $\delta=86=\tilde\delta+1$ and
$\mu_e(\un{r}) = 166 =\frac{1}{2}\cdot \mu_{\tilde e}(\tilde{\un{r}})$,
where again we remark that $\un{r}$ is the minimizing sequence 
associated with $\delta$, and $\tilde{\un{r}}$ is the maximizing
sequence associated with $\tilde\delta$.
Thus $\un{r}$ and $\tilde{\un{r}}$ 
give rise to groups $G$ and $\tilde G$, respectively,
by (\ref{genuscor}) having the same spectrum, precisely one of them 
fulfilling the `$e'=e$' property, and having the same exponent $16$.

\section{\bf Examples: Small rank}\label{smallranksec}

In the remaining two sections, in order to show that the 
combinatorial machinery developed in Section \ref{transsec}
actually is an efficient technique to find $\mu_0(G)$, and in
suitable cases even all of $\textrm{sp}_0(G)$, we explicitly 
work out some `small' examples. Moreover, we show that Talu's
Conjecture (including the case $p=2$) holds within the various
classes of $p$-groups considered.
In this section, now, we deal with the abelian $p$-groups of
minimum genus at most $1$, and those of rank at most $2$, where
in particular we are interested in finding the smallest positive
reduced genus of these groups. 

\subsection{Non-positive reduced minimum genus}\label{nonposex}
We determine the non-trivial abelian $p$-groups $G$ such that
$\mu(G)\in\{0,1\}$, that is $\mu_0(G)\in\{-1,-\frac{1}{2},0\}$.

We have $\mu_i\leq 0$, for $i\in\mathcal I(G)$, if and only if 
$$ {\frac {p-1}{2}} \cdot \wp(s_1,\dotsc, s_i) 
\leq p^i - \lfloor {\frac {s_{i+1}}{2}} \rfloor .$$
From $s_1\geq\cdots\geq s_i\geq 2\cdot\lfloor {\frac {s_{i+1}}{2}} \rfloor+2$
we get
$$(\lfloor {\frac {s_{i+1}}{2}} \rfloor+1)\cdot(p^i-1)\leq 
\frac {p-1}{2} \cdot \wp(s_1,\dotsc, s_i) ,$$
hence assuming $\mu_i\leq 0$ yields 
$$ (\lfloor {\frac {s_{i+1}}{2}} \rfloor+1)\cdot(p^i-1)
\leq p^i - \lfloor {\frac {s_{i+1}}{2}} \rfloor ,$$
that is $\lfloor {\frac {s_{i+1}}{2}} \rfloor \cdot p^i \leq 1$, 
a contradiction for $1\leq i\leq e-1$. 
We consider the remaining cases:
For $i=0$ we get $\mu_0\leq 0$ if and only if 
$\lfloor {\frac {s_1}{2}} \rfloor\leq 1$,
or equivalently $2\leq s_1 \leq 3$, yielding
the cases as indicated in the first table in Table \ref{nonpostbl}, 
where $1\leq e'<e$.
For $i=e$ we get $\mu_e\leq 0$ if and only if 
$\frac {p-1}{2} \cdot \wp(s_1,\dotsc,s_e)\leq p^e$, hence, 
since $s_1\geq\cdots\geq s_e\geq 2$ implies
$p^e-1=\frac {p-1}{2} \cdot \wp(2,\dotsc,2)
\leq \frac {p-1}{2} \cdot \wp(s_1,\dotsc,s_e)$,
we get the cases indicated in the second table in Table \ref{nonpostbl}.

In conclusion, we have $\mu_0(G)<0$, that is $\mu(G)=0$, if and only if 
$$ G\cong \Z_{p^e}\quad\textrm{or}\quad G\cong \Z_2^2 ,$$
and $\mu_0(G)=0$, that is $\mu(G)=1$, if and only if
$$ G\cong\Z_{p^{e'}}\oplus\Z_{p^e}\textrm{ for }e'<e,
\quad\textrm{or}\quad 
G\cong \Z_{p^e}^2 \textrm{ for } p^e\neq 2,
\quad\textrm{or}\quad 
G\cong \Z_2^3 .$$
This also yields all abelian $p$-groups having a genus $g\leq 1$.
Note that the explicit cases for $p=2$ and $p=3$ are precisely the
non-cyclic abelian groups of order at most $9$, which are
treated as exceptional cases in \cite[Thm.4]{mac}.

These results compare to the well-known description of 
finite group actions on compact Riemann surfaces of genus $g\leq 1$, 
see \cite[App.]{sah} or \cite[Sect.6.7]{dou}, as follows:
The cases of $\mu_e<0$ are precisely the abelian $p$-groups
amongst the groups with signature of positive curvature, and belong
to branched self-coverings of the Riemann sphere.
The cases of $\mu_0=0$ and $\mu_e=0$ are precisely the abelian $p$-groups
being smooth epimorphic images of the groups with 
finite signature of zero curvature,
the former belong to unramified coverings of surfaces of 
genus $1$, the latter belong to branched coverings of 
the Riemann sphere by surfaces of genus $1$.

\begin{table}[t]\caption{Non-positive reduced minimum genus.}\label{nonpostbl}
\vspace*{-1.5em}

$$ \begin{array}{|r|l||r|r|}
\hline
 & G & \mu_0 & p\textrm{-datum} \\ 
\hline
\hline
s_1=s_e=2 & \Z_{p^e} & 0 & (0,\dotsc,0;1) \rule{0em}{1.2em} \\
s_1=s_e=3 & \Z_{p^e}^2 & 0 & (0,\dotsc,0;1) \rule{0em}{1.2em} \\
s_1=3>s_e=2 & \Z_{p^{e'}}\oplus\Z_{p^e} & 0 & (0,\dotsc,0;1) 
\rule{0em}{1.2em} \\
\hline
\end{array} $$

$$ \begin{array}{|lr|l||r|r|}
\hline
 & & G & \mu_e & p\textrm{-datum} \\ 
\hline
\hline
 & s_1=s_e=2 & \Z_{p^e} & -1 & (0,\dotsc,0,2;0) \rule{0em}{1.2em} \\
\hline
p=3,\, e=1,& s_1=3 & \Z_3^2 & 0 & (3;0) \rule{0em}{1.2em} \\
\hline
p=2,\, e=1,& s_1=3 & \Z_2^2 & -\frac{1}{2} & (3;0) \rule{0em}{1.2em} \\
p=2,\, e=1,& s_1=4 & \Z_2^3 & 0 & (4;0) \rule{0em}{1.2em} \\
p=2,\, e=2,& s_1=3>s_2=2 & \Z_2\oplus\Z_4 & 0 & (1,2;0) \rule{0em}{1.2em} \\
\hline
\end{array} $$

\vspace*{0.5em}
\hrulefill
\end{table}

\subsection{Groups of rank at most $2$}\label{rank2}
The cases occurring in (\ref{nonposex}) consist 
of all the non-trivial abelian $p$-groups of rank at most $2$,
and the group $G\cong\Z_2^3$. The latter being covered by (\ref{genuscor}), 
we proceed to consider the former in more detail, 
and determine their smallest positive reduced genus $\mu_0^+(G)$, 
and thus their smallest genus $\mu^+(G)\geq 2$. 
The results are collected in Table \ref{rank2tbl}, grouped into 
three infinite series, where $1\leq e'<e$, and finitely many exceptional
cases for $p=2$ and $p=3$. 
The proofs for the cyclic cases and the cases of rank $2$ 
are given in (\ref{cycex}) and (\ref{rk2ex}), respectively;
the cases with $e\leq 2$ will reappear in Section \ref{smallexpsec}.

For the cyclic cases we recover the results in
\cite{har} and \cite[Prop.3.3]{kul}. 
Moreover, we conclude that a cyclic $p$-group is uniquely
determined by its smallest genus $\mu^+(G)\geq 2$, with the single
exception of the groups $\{\Z_2,\Z_4,\Z_8\}$, which indeed have the 
same spectrum $\N_0$.
In particular, Talu's Conjecture (including the case $p=2$)
holds within the class of cyclic $p$-groups.

For the cases of rank $2$ the sharp bound derived here improves
the general bound given in \cite[Prop.3.4]{kul};
and for the cases of cyclic deficiency $\delta=1$, where $p$ is odd, we 
recover the relevant part of \cite[Thm.5.4]{mta} and \cite[Cor.5.5]{mta}.
Moreover, we conclude that an abelian $p$-groups of rank $2$ 
is uniquely determined by its smallest genus $\mu^+(G)\geq 2$, with the 
single exception of the groups $\{\Z_2\oplus\Z_4,\Z_2\oplus\Z_8,\Z_4^2\}$;
it will be shown in (\ref{z4p2spec}) that $\Z_2\oplus\Z_4$ and
$\Z_2\oplus\Z_8$ indeed have the same spectrum $1+2\N_0$, 
which differs from the one of $\Z_4^2$. In particular, Talu's Conjecture 
(including the case $p=2$) holds within the class of abelian
$p$-groups of rank $2$.

\begin{table}[t]\caption{Groups of rank at most $2$.}\label{rank2tbl}
\vspace*{-1.5em}

$$ \begin{array}{|l|l||r|r|}
\hline
G & \textrm{condition} & \mu_0^+(G) & \mu^+(G) \\
\hline
\hline
\Z_{p^e} & p^e\neq 2,3,4 & 
\frac {1}{2}\cdot(p^e-p^{e-1})-1 & 
\frac {1}{2}\cdot(p^e-p^{e-1}) \rule{0em}{1.2em} \\ 
\Z_{p^{e'}}\oplus\Z_{p^e} & (p^{e'},p^e)\neq (2,4) &
\frac{1}{2}\cdot (p^e-p^{e-e'}) -1 &
\frac{1}{2}\cdot p^{e'}\cdot(p^e-p^{e-e'}-2)+1 \rule{0em}{1.2em} \\
\Z_{p^e}^2 & p^e\neq 2,3 & \frac{1}{2}\cdot (p^e-3) & 
\frac{1}{2}\cdot p^e\cdot(p^e-3)+1 \rule{0em}{1.2em} \\
\hline
\end{array} $$ 

$$ \begin{array}{|l||r|r|}
\hline
G & \mu_0^+(G) & \mu^+(G) \\
\hline
\hline
\Z_2 & 1 & 2 \rule{0em}{1.2em} \\
\Z_4 & 1 & 2 \rule{0em}{1.2em} \\
\Z_2^2 & \frac{1}{2} & 2 \rule{0em}{1.2em} \\
\Z_2\oplus\Z_4 & 1 & 3 \rule{0em}{1.2em} \\
\hline
\end{array}
\quad \quad \quad
\begin{array}{|l||r|r|}
\hline
G & \mu_0^+(G) & \mu^+(G) \\
\hline
\hline
\Z_3 & 1 & 2 \rule{0em}{1.2em} \\
\Z_3^2 & 1 & 4 \rule{0em}{1.2em} \\
\hline
\end{array} $$ 

\vspace*{0.5em}
\hrulefill
\end{table}

\subsection{Cyclic groups}\label{cycex}
Let $G\cong\Z_{p^e}$, that is $(s_1,\dotsc,s_e)=(2,\dotsc,2)$;
hence we have $\mathcal I(G)=\{0,e\}$.
By (\ref{nonposex}), we have $\min~\gamma (A_0)=\mu_0=0$ and 
$\min~\gamma (A_e)=\mu_e=-1$, hence both $g=0$ are $g=1$ are genera of $G$. 

We proceed to determine $\mu_0^+(G)$: We have 
$$ A_0=\{(2a,\dotsc,2a)~:~a\geq 1\} ,$$
and hence $\gamma(2a,\dotsc,2a)=(a-1)\cdot p^e$ yields
$$ \min~(\gamma (A_0)\setminus\{0\})=p^e .$$
For $1\leq i\leq e-1$, using the notation of (\ref{minrem}), 
we have $i'=i''=0$ and $\epsilon_i=2$, thus we have $\mu_i=0$ and 
$$ \min~\gamma (A_i)=p^e-p^{e-i}\geq p^e-p^{e-1}=\min~\gamma (A_1) .$$
Moreover, for $p=2$ we have $e'=0$ and
$\min~\gamma (A'_i)=\min~\gamma (A_i)$.
Now let $i=e$:

{\bf i)}
Let first $p$ be odd. Then we have 
$$ A_e=\{(a_1,\dotsc,a_e,2a)~:~a_1\geq\cdots\geq a_e\geq 2(a+1)\} ,$$
hence comparing
$\gamma(a_1,\dotsc,a_e,2a)=-p^e + a +
{\frac {p-1}{2}} \cdot \wp(a_1, \dotsc, a_e)$
with $\gamma(2,\dotsc,2,0)=\mu_e=-1$ yields
$$ \min~(\gamma (A_e)\setminus\{-1\})=
\frac{1}{2}\cdot p^{e-1}\cdot(p-1)-1\geq 0, $$
being attained precisely for $(3,2,\dotsc,2,0)$.
We have $p^{e-1}\cdot (p-1)=2$ if and only if $p=3$ and $e=1$. 
Thus, if $G\not\cong\Z_3$, then we have 
$\mu_0^+(G)=\frac{1}{2}\cdot p^{e-1}\cdot(p-1)-1$.
The case $G\cong\Z_3$ is covered by (\ref{genuscor}).


{\bf ii)}
Let now $p=2$. We have 
$$ A'_e=\{(a_1,\dotsc,a_e,2a)~:~ 
a_1\geq\cdots\geq a_e\geq 2(a+1),\, a_e\textrm{ even}\} .$$
We first assume that $e\geq 3$. Comparing
$\gamma(a_1,\dotsc,a_e,2a)=-2^e + a +
{\frac {1}{2}} \cdot \wp(a_1, \dotsc, a_e)$ 
with $\gamma(2,\dotsc,2,0)=\mu_e=-1$ we get
$$ \min~(\gamma (A'_e)\setminus\{-1\})=2^{e-2}-1>0 ,$$
being attained precisely for $(3,2,\dotsc,2,0)$.
Hence we conclude $\mu_0^+(G)=2^{e-2}-1$.

In particular, for $e=3$, that is $G\cong\Z_8$, we have 
$\gamma(a_1,2,2,0)=2a_1-5$ for $a_1\geq 2$,
and $\gamma(a_1,3,2,0)=2a_1-4$ for $a_1\geq 3$,
implying that $\gamma (A'_3)=\{-1\}\cup\N$.
Hence the reduced spectrum equals $\textrm{sp}_0(\Z_8)=\{-1\}\cup\N_0$, 
yielding the spectrum $\textrm{sp}(\Z_8)=\N_0$;
hence in particular we recover a special case of \cite[Cor.6.3]{kma}.

The case $G\cong\Z_2$ being covered by (\ref{genuscor}), it remains
to consider $G\cong\Z_4$: We have
$$ A'_2=\{(a_1,a_2,2a)~:~ a_1\geq a_2\geq 2(a+1),\, 
  a_2\textrm{ even}\} $$ 
and $\gamma(a_1,a_2,2a)=-4 + a + a_1+\frac{a_2}{2}$.
This yields $\min~(\gamma (A'_2)\setminus\{-1\})=0$,
being attained precisely for $(3,2,0)$, and
$\min~(\gamma (A'_2)\setminus\{-1,0\})=1$,
being attained precisely for $(4,2,0)$.
Thus we have $\mu_0^+(\Z_4)=1$.
From $\gamma(a_1,2,0)=a_1-3$, for $a_1\geq 2$, 
we conclude that $\gamma (A'_2)=\{-1\}\cup\N_0$,
thus the reduced spectrum is $\textrm{sp}_0(\Z_4)=\{-1\}\cup\N_0$,
yielding the spectrum $\textrm{sp}(\Z_4)=\N_0$.


\subsection{Groups of rank $2$}\label{rk2ex}
Let $G\cong\Z_{p^{e'}}\oplus\Z_{p^e}$ for some $1\leq e'\leq e$, 
where for $e'=e$ we get $G\cong\Z_{p^e}^2$; hence 
$(s_1,\dotsc,s_{e'},s_{e'+1},\dotsc,s_e)=(3,\dotsc,3,2,\dotsc,2)$
and $\mathcal I(G)=\{0,e\}$.
By (\ref{nonposex}), we have $\min~\gamma (A_0)=\mu_0=0$, while
$\min~\gamma (A_e)=\mu_e<0$ only for $G\cong\Z_2^2$.
Hence $g=1$ is a genus, while $g=0$ is so if and only if
$G\cong\Z_2^2$. 

We proceed to determine $\mu_0^+(G)$:
We have 
$$ A_0=\{(2a,\dotsc,2a)~:~a\geq 1\} ,$$
and hence $\gamma(2a,\dotsc,2a)=(a-1)\cdot p^e$ yields
$$ \min~(\gamma (A_0)\setminus\{0\})=p^e .$$
Let $1\leq i\leq e-1$. Using the notation of (\ref{minrem}),
for $1\leq i\leq e'$ we have 
$$ \mu_i=-p^e+p^{e-i}\cdot(1+\frac{3}{2}\cdot(p^i-1))
=\frac{1}{2}\cdot p^{e-i}\cdot(p^i-1) ,$$
hence from $i'=0$ and $\epsilon_i=1$ we get
$$ \min~\gamma (A_i)=\mu_i+\frac{1}{2}\cdot p^{e-i}\cdot(p^i-1) 
=p^{e-i}\cdot(p^i-1) .$$
For $e'<i\leq e-1$ we have
$$ \mu_i=-p^e+p^{e-i}\cdot(p^i+\frac{1}{2}\cdot p^{i-e'}\cdot (p^{e'}-1))
=\frac{1}{2}\cdot p^{e-e'}\cdot(p^{e'}-1) ,$$
hence from $i'=e'$ and $i''=0$, as well as $\epsilon_i=2$, we get
$$ \min~\gamma (A_i)=\mu_i+\frac{1}{2}\cdot p^{e-e'}\cdot (p^{e'}+1)-p^{e-i}
=p^{e-i}\cdot(p^i-1) .$$
Thus for all $1\leq i\leq e-1$ we have 
$$ \min~\gamma (A_i)= p^e-p^{e-i} \geq p^e-p^{e-1} = \min~\gamma (A_1) .$$
Moreover, for $p=2$ and $e'<i\leq e-1$ 
we have $\min~\gamma (A'_i)=\min~\gamma (A_i)$.

Hence let $i=e$. We have 
$$ \min~\gamma (A_e)=\mu_e=
-1+\frac{1}{2}\cdot p^{e-e'}\cdot(p^{e'}-1) ,$$
where $\mu_e\leq 0$ if and only if $p^{e-e'}\cdot(p^{e'}-1)\leq 2$,
which holds if and only if $e'=1$ and $p^e\in\{2,3,4\}$.
Hence for $p^e>4$, or $p^e=4$ and $e'=e$, we have $\mu_e>0$.

Assume that $p^e-p^{e-1}<\mu_e=-1+\frac{1}{2}\cdot(p^e-p^{e-e'})$, 
then we have $p^e\cdot(1-\frac{2}{p}+\frac{1}{p^{e'}})<-2$,
implying that $1-\frac{2}{p}+\frac{1}{p^{e'}}<0$,
or equivalently $\frac{2}{p}-\frac{1}{p^{e'}}>1$, a contradiction.

Thus for $1\leq e'<e$ and $(p^{e'},p^e)\neq (2,4)$ we conclude that
$$ \mu_0^+(\Z_{p^{e'}}\oplus\Z_{p^e})
=-1+\frac{1}{2}\cdot p^{e-e'}\cdot(p^{e'}-1) ,$$
and for $p^e\geq 4$ we have 
$$ \mu_0^+(\Z_{p^e}^2)=\frac{1}{2}\cdot (p^e-3).$$

The exceptional cases $G\cong\Z_2\oplus\Z_4$ and 
$G\cong\Z_2^2$ and $G\cong\Z_3^2$ are covered by (\ref{genuscor}).






\subsection{Small $2$-groups}\label{z4p2spec}
As it turns out, the above results already cover all non-trivial 
abelian $2$-groups of order at most $8$.
We observe that in all of these cases there is no spectral gap. 
But this is different for the groups of order $16$, where we have the 
following cases not covered by (\ref{genuscor}):

{\bf i)}
Let $G\cong\Z_4^2$, hence $e'=e=2$, that is $(s_1,s_2)=(3,3)$.
We have seen in (\ref{rk2ex}) that 
$\gamma(A_0)=4\N_0$ and $\min~\gamma (A_1)=2$.
Moreover, we have $\min~\gamma (A_2)=\frac{1}{2}$, where
$$ A_2=\{(a_1,a_2,2a) ~:~ a_1\geq a_2\geq\textrm{max}\{3,2(a+1)\}\}$$
and $\gamma(a_1,a_2,2a)=\frac{1}{2}\cdot(-8+2a+2a_1+a_2)$.
Writing $m\in\N$ as
$$ m=\left\{\begin{array}{ll}
-8+2\cdot\frac{m+5}{2}+3, & \textrm{if }m\textrm{  odd}, \\ 
-8+2\cdot\frac{m+4}{2}+4, & \textrm{if }m\textrm{  even}, \rule{0em}{1.2em} \\ 
\end{array}\right. $$
shows that any $m\in\N\setminus\{2\}$ is of the form $m=-8+2a_1+a_2$ for
some $a_1\geq a_2\geq 3$, while $2$ is not of the form $-8+2a+2a_1+a_2$
for any $(a_1,a_2,2a)\in A_2$. Thus we have
$\gamma(A_2)=(\frac{1}{2}\N)\setminus\{1\}$, hence we conclude that 
$$ \textrm{sp}_0(\Z_4^2)=(\frac{1}{2}\N_0)\setminus\{1\}
\quad\textrm{and}\quad
\textrm{sp}(\Z_4^2)=(1+2\N_0)\setminus\{5\} .$$

{\bf ii)}
Let $G\cong\Z_2\oplus\Z_8$, hence we have $e'=1$ and $e=3$,
that is $(s_1,s_2,s_3)=(3,2,2)$.
We have seen in (\ref{rk2ex}) that $-1$ is not a reduced genus,
and that $\gamma(A_0)=8\N_0$. Moreover, we have $\min~\gamma (A'_3)=1$, where  
$$ A'_3=\{(a_1,a_2,a_3,2a)~:~ 
a_1\geq\textrm{max}\{3,a_2\},\,
a_2\geq a_3\geq 2(a+1),\, a_3\textrm{ even}\} ,$$
and $\gamma(a_1,a_2,a_3,2a) =-8+a+2a_1+a_2+\frac{a_3}{2}$.
Writing $m\in\N$ as
$$ m=\left\{\begin{array}{ll}
-7+2\cdot\frac{m+5}{2}+2, & \textrm{if }m\textrm{  odd}, \\ 
-7+2\cdot\frac{m+4}{2}+3, & \textrm{if }m\textrm{  even}, \rule{0em}{1.2em} \\ 
\end{array}\right. $$
shows that $m=-8+2a_1+a_2+\frac{2}{2}$ for
some $a_1\geq a_2\geq 2$ such that $a_1\geq 3$. 
Thus we have $\gamma(A'_3)=\N$, hence we conclude that
$$ \textrm{sp}_0(\Z_2\oplus\Z_8)=\N_0
\quad\textrm{and}\quad
\textrm{sp}(\Z_2\oplus\Z_8)=1+2\N_0 .$$

{\bf iii)}
Let $G\cong\Z_{16}$, hence we have $e'=0$ and $e=4$,
that is $(s_1,s_2,s_3,s_4)=(2,2,2,2)$.
We have seen in (\ref{cycex}) that $\gamma(A_0)=16\N_0$,
and $\min~\gamma (A'_i)=\min~\gamma (A_i)=16-2^{4-i}$ for $1\leq i\leq 3$.
Moreover, we have $\min~\gamma (A'_4)=-1$, where
$$ A'_4=\{(a_1,a_2,a_3,a_4,2a)~:~ 
a_1\geq a_2\geq a_3\geq a_4\geq 2(a+1),\, a_4\textrm{ even}\} ,$$
and $\gamma(a_1,a_2,a_3,a_4,2a)=-16+a+4a_1+2a_2+a_3+\frac{a_4}{2}$.
Writing $m\in\{-1\}\cup\N_0$ as
$$ m=\left\{\begin{array}{ll}
-15+4\cdot\frac{m+9}{4}+2\cdot 2+2, & \textrm{if }m\equiv 3\pmod{4}, \\ 
-15+4\cdot\frac{m+7}{4}+2\cdot 3+2, & \textrm{if }m\equiv 1\pmod{4}, 
                                      \rule{0em}{1.2em} \\
-15+4\cdot\frac{m+6}{4}+2\cdot 3+3, & \textrm{if }m\equiv 2\pmod{4},
                                      \rule{0em}{1.2em} \\
-15+4\cdot\frac{m+4}{4}+2\cdot 4+3, & \textrm{if }m\equiv 0\pmod{4},
                                      \rule{0em}{1.2em} \\ 
\end{array}\right. $$
shows that any $m\in(\{-1\}\cup\N_0)\setminus\{0,1,2,4,8\}$ 
is of the form $m=-16+4a_1+2a_2+a_3+\frac{2}{2}$ for
some $a_1\geq a_2\geq a_3\geq 2$, 
while none of $\{0,1,2,4,8\}$ is of the form 
$-16+a+4a_1+2a_2+a_3+\frac{a_4}{2}$
for any $(a_1,a_2,a_3,a_4,2a)\in A'_4$. 
Thus we have $\gamma(A'_4)=(\{-1\}\cup\N_0)\setminus\{0,1,2,4,8\}$,
hence we conclude that
$$ \textrm{sp}_0(\Z_{16})=(\{-1\}\cup\N_0)\setminus\{1,2,4\}
\quad\textrm{and}\quad
\textrm{sp}(\Z_{16})=\N_0\setminus\{2,3,5\} ;$$
hence in particular we recover a special case of \cite[Cor.6.3]{kma}.
\QED

\bigskip 
For completeness, the remaining cases 
are dealt with using (\ref{genuscor}), and we get
$$ \textrm{sp}_0(\Z_2^4)=\frac{1}{2}\N 
\quad\textrm{and}\quad
\textrm{sp}_0(\Z_2^2\oplus\Z_4)=\N,  
\quad\textrm{hence}\quad
\textrm{sp}(\Z_2^4)=\textrm{sp}(\Z_2^2\oplus\Z_4)=5+4\N_0 .$$

Collecting the results for all non-trivial abelian $2$-groups of order
at most $16$ yields
$$ \textrm{sp}(\Z_2)=\textrm{sp}(\Z_4)=\textrm{sp}(\Z_2^2)
  =\textrm{sp}(\Z_8)=\N_0 $$
and
$$ \textrm{sp}(\Z_2\oplus\Z_4)=\textrm{sp}(\Z_2^3)
  =\textrm{sp}(\Z_2\oplus\Z_8)=1+2\N_0 .$$
Thus these provide examples of $2$-groups having the
same spectrum, where neither the order, the exponent,
the cyclic deficiency nor the `$e'<e$' property coincide.

\subsection{Small $3$-groups}
By the results above, and (\ref{genuscor}), we have
$$ \textrm{sp}(\Z_3)=\N_0 \quad\textrm{and}\quad
\textrm{sp}(\Z_3^2)=1+3\N_0 \quad\textrm{and}\quad
\textrm{sp}(\Z_3^3)=10+9\N_0 .$$

We again observe that in all of these cases there is no spectral gap,
but this picture already changes for the next $3$-groups springing to 
mind, as soon as we avoid the realm of (\ref{genuscor}). We present a 
couple of examples, showing that going over to reduced spectra tends
to unify and straighten out the computations necessary:

{\bf i)}
Let $G\cong\Z_9$, that is we have $e=2$ and $(s_1,s_2)=(2,2)$.
We have seen in (\ref{cycex}) that 
$\gamma(A_0)=9\N_0$ and $\min~\gamma (A_1)=6$.
Moreover, we have $\min~\gamma (A_2)=-1$, where
$$ A_2=\{(a_1,a_2,2a) ~:~ a_1\geq a_2\geq 2(a+1)\} $$
and $\gamma(a_1,a_2,2a)=-9+a+3a_1+a_2$.
Writing $m\in\N_0$ as
$$ m=\left\{\begin{array}{ll}
-9+3\cdot\frac{m+7}{3}+2, & \textrm{if }m\equiv 2\pmod{3}, \\ 
-9+3\cdot\frac{m+6}{3}+3, & \textrm{if }m\equiv 0\pmod{3}, \rule{0em}{1.2em} \\ 
-9+3\cdot\frac{m+5}{3}+4, & \textrm{if }m\equiv 1\pmod{3}, \rule{0em}{1.2em} \\ 
\end{array}\right. $$
shows that any $m\in\N\setminus\{1,4\}$ can be written as 
$m=-9+3a_1+a_2$ for some $a_1\geq a_2\geq 2$, 
while none of $\{0,1,4\}$ is of the form $-9+a+3a_1+a_2$ 
for any $(a_1,a_2,2a)\in A_2$. Thus we have
$\gamma(A_2)=\{-1\}\cup(\N\setminus\{1,4\})$,
hence we conclude 
$$ \textrm{sp}_0(\Z_9)=(\{-1\}\cup\N_0)\setminus\{1,4\}
\quad\textrm{and}\quad
\textrm{sp}(\Z_9)=\N_0\setminus\{2,5\} ;$$
hence in particular we recover a special case of \cite[Cor.5.3]{kma}.

{\bf ii)}
We determine the spectrum of $G\cong\Z_3\oplus\Z_9$,
thus recovering \cite[Cor.5.5]{mta}:
We have $e'=1$ and $e=2$, that is $(s_1,s_2)=(3,2)$,
and thus we have seen in (\ref{rk2ex}) that 
$\gamma(A_0)=9\N_0$ and $\min~\gamma (A_1)=6$.
Moreover, we have $\min~\gamma (A_2)=2$, where
$$ A_2=\{(a_1,a_2,2a) ~:~ a_1\geq\textrm{max}\{3,a_2\},\,a_2\geq 2(a+1)\}$$
and $\gamma(a_1,a_2,2a)=-9+a+3a_1+a_2$. As above, writing $m\in\N$ as
$$ m=\left\{\begin{array}{ll}
-9+3\cdot\frac{m+7}{3}+2, & \textrm{if }m\equiv 2\pmod{3}, \\ 
-9+3\cdot\frac{m+6}{3}+3, & \textrm{if }m\equiv 0\pmod{3}, \rule{0em}{1.2em} \\ 
-9+3\cdot\frac{m+5}{3}+4, & \textrm{if }m\equiv 1\pmod{3}, \rule{0em}{1.2em} \\ 
\end{array}\right. $$
shows that any 
$m\in\N\setminus\{1,4\}$ can be written as $m=-9+3a_1+a_2$ for
some $a_1\geq\textrm{max}\{3,a_2\}$ and $a_2\geq 2$, 
while none of $\{1,4\}$ is of the form $-9+a+3a_1+a_2$
for any $(a_1,a_2,2a)\in A_2$.
Thus we have $\gamma(A_2)=\N\setminus\{1,4\}$. Hence we conclude that 
$$\textrm{sp}_0(\Z_3\oplus\Z_9)=\N_0\setminus\{1,4\}
\quad\textrm{and}\quad
\textrm{sp}(\Z_3\oplus\Z_9)=(1+3\N_0)\setminus\{4,13\} .$$

\section{\bf Examples: Small exponents}\label{smallexpsec}

In this section we consider abelian $p$-groups of exponent
at most $p^2$. In particular, we ask ourselves whether the description
of the reduced minimum genus in terms of the defining invariants of the
group in question lends itself to a `generic' description.

\subsection{Elementary abelian groups}\label{e1ex}
Let $G\cong\Z_p^r$ be an elementary abelian $p$-group, that is $e=1$,
and let $s:=r+1\geq 2$. We have and $\mathcal I(G)=\{0,1\}$, where
(\ref{minthm}) says that $0\in\mathcal I(G)$ can be
ignored whenever $s$ is even. Still, we have
$$ {\mathrm {min}} ~\gamma (A_0)=\mu_0 =
\left\{\begin{array}{ll}
\frac{ps}{2}-p,& \textrm{if }s\textrm{ even}, \\
\frac{ps}{2}-\frac{3p}{2},& \textrm{if }s\textrm{ odd}, \rule{0em}{1.2em} \\
\end{array}\right. $$
and 
$$ {\mathrm {min}} ~\gamma (A_1)=\mu_1=\frac{ps}{2}-\frac{s}{2}-p .$$
Thus we have $\mu_0<\mu_1$ if and only if $s$ is odd and $s<p$,
with equality if and only if $s=p$ is odd.
Hence we get $\mu_0(G)=\mu_0$ if $s$ is odd and $s<p$,
otherwise we have $\mu_0(G)=\mu_1$.
In particular, for $p$ odd we thus recover, and at the same time 
correct \cite[Sect.7, Rem.]{mta}, where $\mu_0(G)$ is erroneously stated 
for $s<p$.

We are tempted to call the cases where $s$ is odd such that $s<p$
the `exceptional' ones, and the remaining the `generic' ones;
then there are only finitely many `exceptional' cases,
which do not occur at all for $p=2$.
In particular, as part of the `generic' region 
we have $\mu_0(G)=\mu_1$ for $s\geq\max\{p-1,2\}$, 
in accordance with (\ref{genusthm}).

{\bf i)}
For $p$ odd, viewing $\mu_0$ and $\mu_1$ as linear functions in $s$, 
with positive slope $\frac{p}{2}$ and $\frac{p-1}{2}$, respectively, 
and since $\mu_0(s+1)-\mu_1(s)=\frac{s}{2}>0$, for $2\leq s<p$ even,
we conclude that the reduced minimum genus $\mu_0(G)$ is strictly increasing
with $s$, and thus the minimum genus $\mu(G)=1+p^{s-2}\cdot\mu_0(G)$ is
as well. 

{\bf ii)}
For $p=2$ we have $\mu_0(G)=\mu_1=\frac{s}{2}-2$ 
for all $s\geq 2$, thus the reduced minimum genus $\mu_0(G)$ 
is strictly increasing with $s$, and hence the minimum genus
$\mu(G)=1+2^{s-2}\cdot\mu_0(G)$, for $s\geq 3$, is as well. 

A few values are given in the first and second table in 
Table \ref{elemtbl}, respectively.
We conclude that $G$ is uniquely determined by its minimum
genus $\mu(G)$, with the single exception of $\{\Z_2,\Z_2^2\}$;
indeed, as we have already noted in (\ref{z4p2spec}), 
the latter have the same spectrum.
Hence in particular Talu's Conjecture (including the
case $p=2$) holds within the class of elementary abelian $p$-groups.

Note that for $p$ odd this would also be a consequence of 
\cite[Cor.7.3]{mta}, but due to the
erroneous \cite[Sect.7, Rem.]{mta} the results \cite[Thm.7.2, Cor.7.3]{mta}
are at stake; only \cite[Cor.7.3(1)]{mta} can be verified 
independently by (\ref{genusthm}).

\begin{table}[t]\caption{Elementary abelian groups}\label{elemtbl}
\vspace*{-1.5em}

$$ \begin{array}{|r||cccccccccc|}
\hline
s & 2 & 3 & 4 & 5 & \!\dotsc\! & p-3 & p-2 & p-1 & p & p+1 \\
\hline
\mu_0(G) & -1 & 0 & p-2 & p & \!\dotsc\! & 
\frac{p(p-6)+3}{2} & \frac{p(p-5)}{2} & \frac{p(p-4)+1}{2} &
\frac{p(p-3)}{2} & \frac{p(p-2)-1}{2} \rule{0em}{1.2em} \\
\hline
\end{array} $$

$$ \begin{array}{|r||rrrrrrr|}
\hline
s & 2 & 3 & 4 & 5 & 6 & 7 & 8 \\
\hline
\mu_0(G) & -1 & -\frac{1}{2} & 0 & \frac{1}{2} & 1 & \frac{3}{2} & 2 
\rule{0em}{1.2em} \\
\mu(G) & 0 & 0 & 1 & 5 & 17 & 49 & 129 \rule{0em}{1.2em} \\
\hline
\end{array} $$

\vspace*{0.5em}
\hrulefill
\end{table}

\subsection{Groups of exponent $p^2$}\label{e2ex}
Let $G\cong\Z_p^{r_1}\oplus\Z_{p^2}^{r_2}$, that is we have $e=2$.
Let $s:=s_1=r_1+r_2+1$ and $t:=s_2=r_2+1$, hence $s\geq t\geq 2$.
Moreover, we have $\{0,2\}\subseteq\mathcal I(G)\subseteq\{0,1,2\}$,
where $1\in\mathcal I(G)$ if and only if 
$s-t\geq 2$, or $s-t=1$ and $t$ is odd; additionally, (\ref{minthm}) 
says that $0\in\mathcal I(G)$ can be ignored whenever $s$ is even.

Still, in order to obtain a complete overview, we explicitly have
$$ {\mathrm {min}} ~\gamma (A_0)=\mu_0 =
\left\{\begin{array}{ll}
\frac{p^2s}{2}-p^2,& \textrm{if }s\textrm{ even}, \\
\frac{p^2s}{2}-\frac{3p^2}{2},& \textrm{if }s\textrm{ odd}, 
\rule{0em}{1.2em} \\
\end{array}\right. $$
and 
$$ {\mathrm {min}} ~\gamma (A_1)=
\left\{\begin{array}{ll}
\frac{p^2s}{2}-\frac{p(s-t)}{2}-p^2 ,& 
\textrm{if }t\textrm{ even},\, s-t\geq 2, \\
\frac{p^2s}{2}-\frac{p(s-t)}{2}-\frac{p}{2}-p^2 ,& 
\textrm{if }t\textrm{ odd},\, s-t\geq 2, \rule{0em}{1.2em} \\
\frac{p^2s}{2}-p-\frac{p^2}{2} ,& 
\textrm{if }t\textrm{ even},\, s-t=1, \rule{0em}{1.2em} \\
\frac{p^2s}{2}-p-p^2 ,& 
\textrm{if }t\textrm{ odd},\, s-t=1, \rule{0em}{1.2em} \\
\frac{p^2s}{2}-p ,& 
\textrm{if }t\textrm{ even},\, s=t, \rule{0em}{1.2em} \\
\frac{p^2s}{2}-p-\frac{p^2}{2} ,& 
\textrm{if }t\textrm{ odd},\, s=t, \rule{0em}{1.2em} \\
\end{array}\right. $$
and
$$ {\mathrm {min}} ~\gamma (A_2)=\mu_2= 
\frac{p^2s}{2}-\frac{p(s-t)}{2}-\frac{t}{2}-p^2 .$$

Thus we have
$\mu_0<\mu_2$ if and only if $s$ is odd and $p(s-t)+t<p^2$,
with equality if and only if $s$ is odd and $p(s-t)+t=p^2$;
and $({\mathrm {min}} ~\gamma (A_1))<\mu_2$ 
if and only if $t$ is odd and $t<\min\{p,s\}$,
with equality if and only if $t=p$ is odd and $t<s$;
and $\mu_0<({\mathrm {min}} ~\gamma (A_1))$ if and only if 
$$ \left\{\begin{array}{l}
s=t\textrm{ even}, \\
s\textrm{ odd},\, t\textrm{ even},\, s-t<p, \rule{0em}{1.2em} \\
s\textrm{ odd},\, t\textrm{ odd},\, s-t< p-1, \rule{0em}{1.2em} \\
\end{array}\right. $$
with equality if and only if $s$ is odd, and $s-t=p$ odd or $s-t=p-1$ even.

In particular, we have equality
$\mu_0=({\mathrm {min}} ~\gamma (A_1))=\mu_2$ throughout if and only if
$t=p$ odd and $s=2p-1$.
Anyway, there are three cases in which $\mu_0(G)$ coincides with either 
of $\mu_0$, $\mu_1$ and $\mu_2$ in turn, where the mutual intersection
of these cases is described by equating the associated $\mu_i$:

{\bf i)}
Let $s$ be odd such that $p(s-t)+t\leq p^2$, thus $\mu_0\leq\mu_2$. 
Moreover, we have $s-t<p$, implying 
$\mu_0\leq({\mathrm {min}} ~\gamma (A_1))$, 
hence we get 
$$\mu_0(G)=\mu_0=\frac{p^2}{2}\cdot (s-3) .$$

{\bf ii)}
Let $t$ be odd such that $t\leq p$, and let $s$ be even or $s-t\geq p-1$.
Then we have $({\mathrm {min}} ~\gamma (A_1))\leq\mu_2$ and 
$({\mathrm {min}} ~\gamma (A_1))\leq\mu_0$, hence we get 
$$ \mu_0(G)=\mu_1=\frac{p^2}{2}\cdot (s-2)-\frac{p}{2}\cdot (s-t+1) .$$

{\bf iii)}
Let $s$ be even or $p(s-t)+t\geq p^2$, and 
let $t$ be even or $t=s$ or $t\geq p$. 
Then we have $\mu_2\leq\mu_0$ and 
$\mu_2\leq({\mathrm {min}} ~\gamma (A_1))$, hence we get 
$$ \mu_0(G)=\mu_2= 
\frac{p^2}{2}\cdot (s-2)-\frac{p}{2}\cdot (s-t) - \frac{t}{2} .$$

Note that case i) consists of finitely many pairs $(s,t)$,
while in case ii) $s$ is unbounded but $t$ is still bounded.
Hence we are again tempted to call these the `exceptional' cases, 
as opposed to the `generic' case iii), where both $s$ and $t$ are unbounded.
In particular, as part of the `generic' region we have $\mu_0(G)=\mu_2$ 
for $t\geq\max\{p-1,2\}$ and $s-t\geq p-1$,
which we will recover as a special case of (\ref{genusthm}).
In particular, for $p=2$ case i) consists of the pairs $(s,t)=(3,3)$
and $(s,t)=(3,2)$, that is $G\cong\Z_4^2$ and $G\cong\Z_2\oplus\Z_4$, 
respectively, case ii) does not occur at all, and all pairs except 
$(s,t)=(3,3)$ belong to case iii).

To further illustrate the idea of distinguishing between `generic' and 
`exceptional' pairs, the various cases for $p=5$ and $2\leq t\leq s\leq 27$ 
are visualized in Table \ref{regiontbl}: 
The cases i), ii) and iii) are depicted by 
`$\ast$', `$\bullet$' and `$\cdot$', respectively,
the intersections `i)$\cap$iii)' and `ii)$\cap$iii)' are indicated by 
`$\times$' and `$\circ$', respectively, 
and `i)$\cap$ii)', consisting of $(s,t)\in\{(7,3),(9,5)\}$, 
is indicated by `$\circledast$' and `$\otimes$', where the latter
icon refers to `i)$\cap$ii)$\cap$iii)', which is $(s,t)=(9,5)$.

The closed interior of the cone emanating from $(s,t)=(8,4)$
indicates the realm of applicability of (\ref{genusthm}); actually, 
this turns out to be the largest cone being contained in the `generic' 
region, saying that in a certain sense this result is best possible, 
at least for the cases considered here.
Moreover, within this cone, the `generic' case iii) refers to 
the case $j=2$ in the notation of (\ref{genusthm}), while the
`exceptional' intersection `ii)$\cap$iii)' refers to $j\leq 1$,
that is the pairs $(s,5)$ such that $s\geq 9$, and finally the intersection 
`i)$\cap$ii)$\cap$iii)' refers to $j=0$, that is $(s,t)=(9,5)$. 

\begin{table}\caption{`Generic' and `exceptional' cases for $p=5$.}
\label{regiontbl}
\begin{center}
\includegraphics[height=125mm]{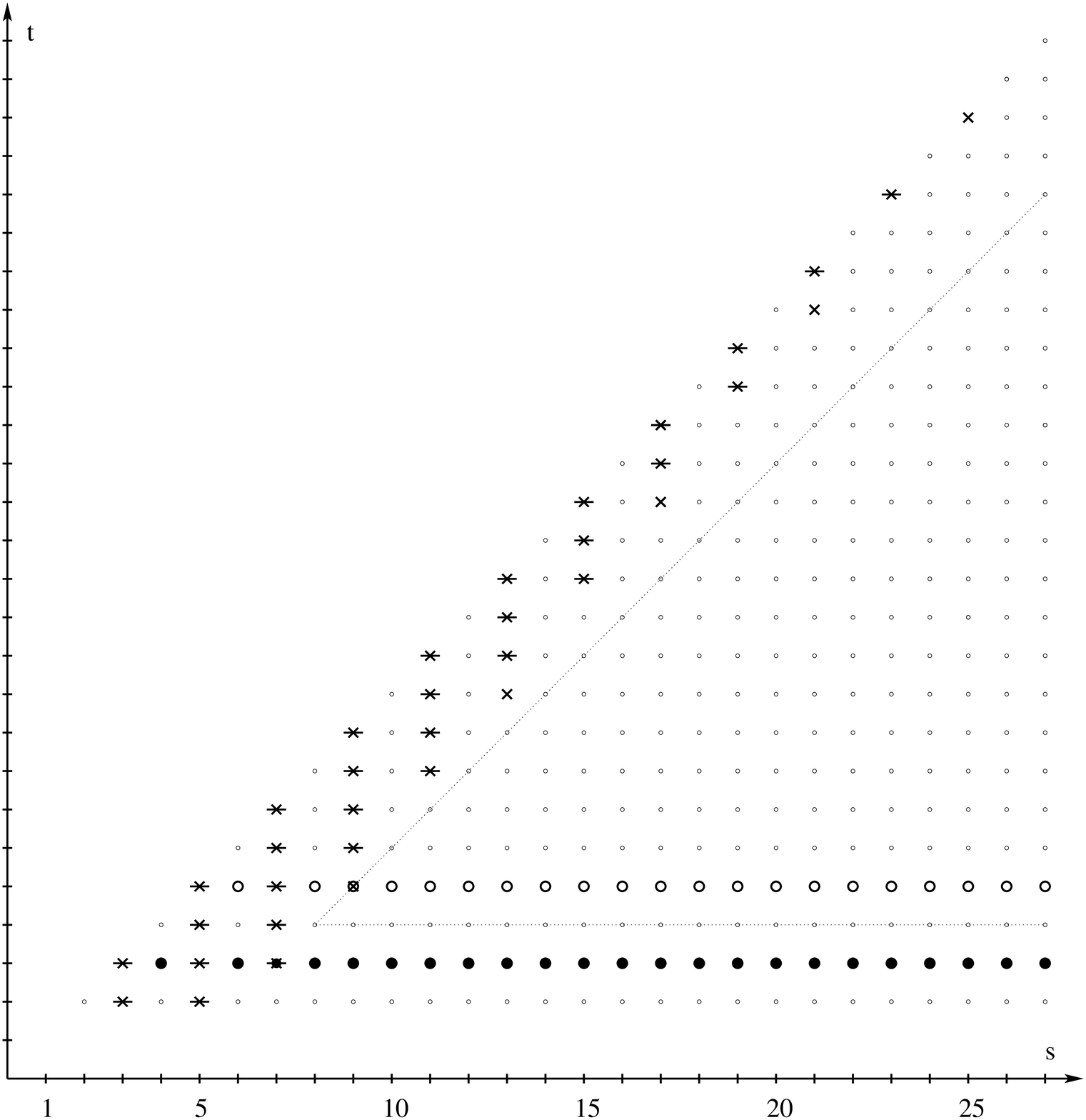}
\end{center}
\hrulefill
\end{table}

\subsection{Recovering groups of exponent $p^2$}\label{e2uniq}
Keeping the notation of (\ref{e2ex}), we show that $G$ is uniquely
determined by its Kulkarni invariant $N=N(G)$ and its minimum genus $\mu(G)$,
with the single exception of the groups $\{\Z_4^2,\Z_2\oplus\Z_4\}$; the 
latter groups can be distinguished by their spectrum, see (\ref{z4p2spec}).
In particular we conclude that Talu's Conjecture (including the
case $p=2$) holds within the class of abelian $p$-groups of exponent $p^2$;
thus for $p$ odd we recover \cite[Thm.3.8]{tal}:

Let first $p$ be odd. The cyclic deficiency $\delta=\delta(G)$ and
the reduced minimum genus $\mu_0(G)$ of $G$ are known from
$\delta=\log_p(N)$ and $\mu_0(G)=\frac{\mu(G)-1}{p^\delta}$. We have 
$\delta=r_1+2r_2-2=s+t-4$, thus we may view 
$\mu_0$ in case (\ref{e2ex}.i), $\mu_1$ in case (\ref{e2ex}.ii), 
and $\mu_2$ in case (\ref{e2ex}.iii)
as linear functions in $s$, depending on the parameter $\delta$:
$$ \begin{array}{lcl}
\mu_0&=&\frac{p^2}{2}\cdot s-\frac{3p^2}{2}, \\
\mu_1&=&(\frac{p^2}{2}-p)\cdot s+\frac{p(\delta+3)}{2}-p^2, 
\rule{0em}{1.2em} \\
\mu_2&=&\frac{(p-1)^2}{2}\cdot s+\frac{(p-1)(\delta+4)}{2}-p^2. 
\rule{0em}{1.2em} \\ 
\end{array}$$
As these functions have positive slope, they are strictly increasing,
hence we look for coincidences across cases:

{\bf i)}
Let first $\mu_1(s,t)=\mu_0(\tilde s,\tilde t)$, where 
$(s,t)$ belongs to case (\ref{e2ex}.ii), 
and $(\tilde s,\tilde t)$ belongs to case (\ref{e2ex}.i).
Then we conclude that $\tilde s=s-\frac{s-t+1}{p}+1$, hence we have
$s-t=kp-1$ for some $k\geq 1$. From this get 
$s=\frac{1}{2}\cdot(\delta+3+kp)$ and $t=\frac{1}{2}\cdot(\delta+5-kp)$,
implying $\tilde s=s-k+1=\frac{1}{2}\cdot(\delta+5-2k+kp)$ and
$\tilde t=\delta+4-\tilde s=\frac{1}{2}\cdot(\delta+3+2k-kp)$.
Thus we get $\tilde s-\tilde t=1+k(p-2)$. Hence $\tilde s-\tilde t\leq p-1$
yields $k=1$, and thus $\tilde s=s$ and $\tilde t=t$.
Note that in this case both $s$ and $t$ are odd such that $t\leq p$ and
$s-t=p-1$, indeed yielding $\mu_1(s,t)=\mu_0(s,t)$. 

{\bf ii)}
Let next $\mu_2(s,t)=\mu_0(\tilde s,\tilde t)$, where 
$(s,t)$ belongs to case (\ref{e2ex}.iii), and $(\tilde s,\tilde t)$
belongs to case (\ref{e2ex}.i).
Then we conclude that $\tilde s=s+\frac{(p-1)t-ps}{p^2}+1$, hence we have
$t=kp$ for some $k\geq 1$. Thus we infer that $p$ divides $k(p-1)-s$, 
hence we get $s=k(p-1)+lp$ for some $l\geq 1$. This yields
$\tilde s=(k+l)(p-1)+1$ and $\tilde t=s-\tilde s+t=kp+l-1$.
Hence we have $p(\tilde s-\tilde t)+\tilde t=l(p-1)^2+2p-1\leq p^2$,
implying $l=1$, thus $\tilde s=s=(k+1)p-k$ and hence $\tilde t=t$.
Note that in this case $s$ is odd, where $s-t=p-k$ and $t=kp\geq p$, 
hence $p(s-t)+t=p^2$, indeed yielding $\mu_2(s,t)=\mu_0(s,t)$. 

{\bf iii)}
Let finally $\mu_2(s,t)=\mu_1(\tilde s,\tilde t)$, where 
$(s,t)$ belongs to case (\ref{e2ex}.iii), and $(\tilde s,\tilde t)$
belongs to case (\ref{e2ex}.ii).
Then we conclude that $(p-1)\tilde s+\tilde t-1=(p-1)s+\frac{p-1}{p}\cdot t$,
hence we have $t=kp$ for some $k\geq 1$, and thus 
$\tilde t-1=(p-1)(s+k-\tilde s)\geq p-1$. This yields
$\tilde s=s+k-1$ and $\tilde t=p$. Hence we get 
$s+kp=s+t=\delta+4=\tilde s+\tilde t=s+k-1+p$,
implying $(k-1)p=k-1$, thus $k=1$, and hence $\tilde s=s$ and $\tilde t=t$.
Note that in this case $t=p$ is odd, and $s$ is even or $s\geq 2p-1$, 
in particular yielding $\mu_2(s,t)=\mu_1(s,t)$.  

This concludes our treatment of the case $p$ odd, hence let now $p=2$. 

{\bf i)}
We first consider case (\ref{e2ex}.iii), where,
using $s+t=\delta-4$ again, we have
$$ \mu_2=s+\frac{t}{2}-4=\frac{s}{2}+\frac{\delta}{2}-2 .$$
We distinguish the cases $t=2$ and $t>2$:
If $t=2$, then we have $\log_2(N)=\delta=s-2=\mu_2+1$, thus 
$$ \mu(G)=\mu_2\cdot 2^\delta+1=(\log_2(N)-1)\cdot N+1 ,$$
while if $t>2$, then we have $\log_2(N)=\delta-1$, thus 
$$ \mu(G)=\mu_2\cdot 2^\delta+1=(\log_2(N)+s-3)\cdot N+1 .$$
Hence we are able to decide in which of these cases we are, 
and to determine $\delta$ and subsequently $s$, in the former
case from $N$, in the latter case from $N$ and $\mu(G)$.

{\bf ii)}
Finally, we consider the pair $(3,3)$, that is $G\cong\Z_4^2$, 
which is the only pair not belonging to case (\ref{e2ex}.iii), 
but just to case (\ref{e2ex}.i):
We have $\mu_0(\Z_4^2)=\mu_0(3,3)=0$, hence its minimum genus equals 
$\mu(\Z_4^2)=1$. For pairs $(s,t)$ belonging to case (\ref{e2ex}.iii), 
the statement $\mu(G)=1$ translates into
$\mu_2(s,t)=0$, that is $s+\frac{t}{2}=4$, being
equivalent to $(s,t)=(3,2)$, that is $G\cong\Z_2\oplus\Z_4$;
note that $(3,2)$ is the other pair belonging to case (\ref{e2ex}.i).
Moreover, for $G\cong\Z_4^2$ we have $\log_2(N)=\delta-1=1$,
and for $G\cong\Z_2\oplus\Z_4$ we also have $\log_2(N)=\delta=1$.
Thus $\{\Z_4^2,\Z_2\oplus\Z_4\}$ are the only groups under 
consideration which cannot be distinguished by $N$ and $\mu(G)$.

\renewcommand\refname


\vspace*{1em}

{J.M.:
Arbeitsgruppe Algebra und Zahlentheorie \\
Bergische Universit\"at Wuppertal \\
Gau{\ss}-Stra{\ss}e 20 \\
D-42119 Wuppertal, Germany} \\
{\sf juergen.mueller@math.uni-wuppertal.de}

\vspace*{1em}

{S.S.:
Department of Mathematics \\
Indian Institute of Science Education and Research Bhopal \\
Indore-Bypass Road, Bhauri \\
Bhopal 462066, India} \\
{\sf sidhu@iiserb.ac.in}

\end{document}